\documentclass[gammmod,a4paper,fleqn%
]{w-artmod}
\usepackage{times,cite,w-thm}
\usepackage{amsmath}
\usepackage{amssymb}
\usepackage{amstext}
\usepackage{amsthm}
\usepackage{a4wide}
\usepackage{color}
\usepackage{url}
\usepackage{hyperref}
\usepackage{stmaryrd}
\usepackage{graphicx}
\usepackage{algorithmic}
\usepackage{algorithm}

\theoremstyle{plain}

\theoremstyle{definition}

\usepackage[]{graphicx}

\newcommand{\calA}{\mathcal{A}}
\newcommand{\calB}{\mathcal{B}}
\newcommand{\calC}{\mathcal{C}}

\newcommand{\calM}{\mathcal{M}}

\newcommand{\calP}{\mathcal{P}}

\newcommand{\calT}{\mathcal{T}}
\newcommand{\calU}{\mathcal{U}}

\newcommand{\calW}{\mathcal{W}}
\newcommand{\calX}{\mathcal{X}}
\newcommand{\calY}{\mathcal{Y}}

\newcommand{\R}{{\mathbb R}}

\DeclareMathOperator{\vect}{vec}

\begin{document}
\DOIsuffix{theDOIsuffix}
\Volume{29}
\Month{01}
\Year{2007}
\pagespan{1}{}
\Receiveddate{XXXX}
\Reviseddate{XXXX}
\Accepteddate{XXXX}
\Dateposted{XXXX}

\title[Low-rank tensor approximation techniques]{A literature survey of low-rank tensor approximation techniques}

\author[L. Grasedyck]{Lars Grasedyck\inst{1,}\footnote{lgr@igpm.rwth-aachen.de}}

\address[\inst{1}]{   Institut f\"ur Geometrie  und Praktische
  Mathematik, RWTH Aachen, Templergraben 55, 52056 Aachen, Germany}
\author[D. Kressner]{Daniel Kressner\inst{2,}\footnote{daniel.kressner@epfl.ch}}
\address[\inst{2}]{ANCHP, MATHICSE, EPF Lausanne, Switzerland}
\author[C. Tobler]{Christine Tobler\inst{2,}\footnote{christine.tobler@epfl.ch}}


\begin{abstract}
During the last years, low-rank tensor approximation has been established as a new tool in
scientific computing to address large-scale linear and multilinear algebra problems, which would be
intractable by classical techniques. This survey attempts to give a literature overview of
current developments in this area, with an emphasis on function-related tensors.
\end{abstract}

\maketitle

\section{Introduction}

This survey is concerned with tensors in the sense of multidimensional arrays.
A general tensor of \emph{order} $d$ and size $n_1\times n_2 \times \cdots \times n_d$ for
integers $n_1,n_2,\ldots,n_d$ will be denoted by
\[
 \calX \in \R^{n_1\times n_2\times \cdots \times n_d}.
\]
An entry of $\calX$ is denoted by $\calX_{i_1,\ldots,i_d}$ where each index $i_\mu \in \{1,\ldots,n_\mu\}$ refers to the $\mu$th mode of the tensor for $\mu = 1,\ldots,d$.
For simplicity, we will assume that $\calX$ has real entries, but it is of course possible to define
complex tensors or, more generally, tensors over arbitrary fields.

A wide variety of applications lead to problems where the data or the desired solution can be represented
by a tensor. In this survey, we will focus on tensors that are induced by the discretization of a multivariate
function; we refer to the survey~\cite{Kolda2009} and to the books~\cite{Kroonenberg2008,Smilde2004} for the treatment of
tensors containing observed data. The simplest way a given multivariate function $f(x_1,x_2,\ldots,x_d)$
on a tensor product domain $\Omega = [0,1]^d$ leads to a tensor is by sampling $f$ on a tensor grid.
In this case, each entry of the tensor contains the function value at the corresponding position in the grid.
The function $f$ itself may, for example, represent the solution to a high-dimensional partial differential equation (PDE).

As the order $d$ increases, the number of entries in $\calX$ increases exponentially for
constant $n = n_1 = \cdots = n_d$. This so called \emph{curse of dimensionality} prevents the
explicit storage of the entries except for very small values of $d$. Even for $n = 2$, storing a tensor
of order $d = 50$ would require 9 petabyte! It is therefore essential to approximate tensors of
higher order in a compressed scheme, for example, a low-rank tensor decomposition. Various such decompositions have been developed, see
Section~\ref{sec:formats}. An important difference to tensors containing observed data, a 
tensor $\calX$ induced by a function is usually not given directly but only as the solution of
some algebraic equation, e.g., a linear system or eigenvalue problem. This requires the development of
solvers for such equations working within the compressed storage scheme. Such algorithms are discussed
in Section~\ref{sec:algorithms}.

The range of applications of low-rank tensor techniques is quickly expanding. For example,
they have been used for addressing:
\begin{itemize}
\item the approximation of parameter-dependent integrals~\cite{Ballani2012b,Khoromskij2011f,Meszmer2012},
multi-dimensional integrals~\cite{Chinnamsetty2012,Hackbusch2007,Khoromskij2011},
and multi-dimensional convolution~\cite{Khoromskij2010b};
\item computational tasks in electronic structure calculations, e.g., based on Hartree-Fock or DFT 
  models~\cite{Benedikt2011,Bertoglio2012,Blesgen2012,Chinnamsetty2007,Chinnamsetty2010a,Chinnamsetty2010,Flad2010,Khoromskij2007a,Khoromskij2008,Khoromskij2009,Khoromskij2009a,Khoromskij2011d,Khoromskaia2010,Khoromskaia2010a,Khoromskaia2011,Khoromskaia2012,Khoromskaia2012c,Khoromskaia2013};
\item the solution of stochastic and parametric PDEs~\cite{Doostan2009,Doostan2012,Espig2013,Espig2013a,Khoromskij2011b,Kressner2011,Khoromskij2010,Matthies2012,Nouy2010};
\item approximation of Green's functions of high-dimensional PDEs~\cite{Biagioni2012,Hackbusch2008a};
\item the solution of the Boltzmann equation~\cite{Ibragimov2009,Khoromskij2007}, and the chemical master / Fokker-Planck equations~\cite{Ammar2006,Ammar2010a,Dolgov2012f,Dolgov2012g,Figueroa2012,Jahnke2008,Hegland2011,Kazeev2013,Mokdad2010}
\item the solution of high-dimensional Schr\"odinger equations~\cite{Chinnamsetty2009,Kressner2011a,Khoromskij2010a,Lubich2008,Meyer1990,Meyer2009};
\item computational tasks in micromagnetics~\cite{Exl2012,Exl2012a};
\item rational approximation problems~\cite{Takaaki2012};
\item computational homogenization~\cite{Chinesta2008,Giraldi2013,Lamari2010};
\item computational finance~\cite{Jondeau2010,Kazeev2012b};
\item the approximation of stationary states of stochastic automata networks~\cite{Buchholz2010};
\item multivariate regression and machine learning~\cite{Beylkin2009}.
\end{itemize}
Note that the list above mainly focuses on techniques that involve tensors; low-rank matrix techniques, such as POD and reduced basis methods, have been applied in an even broader  setting.

A word of caution is appropriate. Even though the field of low-rank tensor approximation is
relatively young, it already appears a daunting task to give proper credit to all developments
in this area. This survey is biased towards work related to the TT and hierarchical Tucker decompositions.
Surveys with a similar scope are the lecture notes by Grasedyck~\cite{Grasedyck2010a}, Khoromskij~\cite{Khoromskij2005,Khoromskij2011,Khoromskij2012},
and Schneider~\cite{Schneider2012}, as well as the monograph by Hackbusch~\cite{Hackbusch2012}.
Less detailed attention may be given to other developments, although we have made an effort to at least touch on all important directions in the area of function-related tensors.

\section{Low-rank tensor decompositions} \label{sec:formats}

As mentioned in the introduction, it will rarely be possible to store
all entries of a higher-order tensor explicitly. Various compression
schemes have been developed to reduce storage requirements. For $d = 2$  all these schemes boil down to the well known reduced singular value
decomposition (SVD) of matrices~\cite{Golub1996}; however, they differ significantly
for tensors of order $d \ge 3$.

\subsection{CP decomposition} \label{sec:cp}

The entries of a rank-one tensor $\calX$ can be written as
\begin{equation} \label{eq:cpentries}
 \calX_{i_1,i_2,\ldots,i_d} = u_{i_1}^{(1)} u_{i_2}^{(2)} \cdots u_{i_d}^{(d)}, \quad 1\le i_\mu \le n_\mu, \quad \mu = 1,\ldots, d.
\end{equation}
By defining the vectors $u^{(\mu)} := \big(u_{1}^{(\mu)},\ldots,u_{n_\mu}^{(\mu)}\big)^T$, a more compact from of this relation is given by
\[
 \vect(\calX) = u^{(d)} \otimes u^{(d-1)} \otimes \cdots \otimes u^{(1)},
\]
where $\otimes$ denotes the usual Kronecker product and $\vect$ stacks the entries of a tensor into a long column vector, such that the indices are in \emph{reverse} lexicographical order.
(Using the vector outer product $\circ$, this relation takes the form
$\calX = u^{(1)}\circ u^{(2)} \circ\cdots \circ u^{(d)}$.)
When $\calX$ represents the discretization of a separable function $f(x_1,x_2,\ldots,x_d) = f_1(x_1)f_2(x_2)\ldots f_d(x_d)$ then $\calX$ is a rank-one tensor with each vector $u^{(\mu)}$ corresponding to a discretization of $f_\mu$.

The \emph{CP (CANDECOMP/PARAFAC) decomposition} is a sum of rank-one tensors:
\begin{equation} \label{eq:def_cp}
\vect(\calX) = 
u^{(d)}_1 \otimes u^{(d-1)}_1 \otimes \cdots \otimes u^{(1)}_1
+ \cdots +
u^{(d)}_R \otimes u^{(d-1)}_R \otimes \cdots \otimes u^{(1)}_R.
\end{equation}
The \emph{tensor rank} of $\calX$ is defined as the minimal $R$ such that $\calX$ has a
CP decomposition with $R$ terms. Note that, in contrast to matrices, the set CP($R$) of
tensors of rank at most $R$ is in general not closed, which renders the problem of finding a best low-rank approximation ill-posed~\cite{Silva2008}. For more properties of the CP decomposition,
we refer to the survey paper~\cite{Kolda2009}.

The CP decomposition requires the storage of $(n_1 + n_2 + \cdots + n_d) R$ entries, which becomes very attractive
for small $R$. To be able to use the CP decomposition for the approximation of function-related tensors, robust and efficient compression techniques are essential.
In particular, the truncation of a rank-$R$ tensor to lower tensor rank is frequently needed.
Nearly all existing algorithms are based on carefully adapting existing optimization algorithms,
see, once again,~\cite{Kolda2009} for an overview of the literature until around 2009. More recent developments
for general tensors include work on
increasing the efficiency and robustness of gradient-based and Newton-like methods~\cite{Acar2011,Espig2008,Espig2012,Kazeev2010,Phan2012b,Phan2013},
modifying and improving ALS (alternating least squares)~\cite{Chen2011,Friedland2012,DeSterck2011,Sterck2012,Rajih2008},
studying the convergence of ALS~\cite{Comon2009,Mohlenkamp2013,Uschmajew2012a}
and reducing the cost of the unfolding operations required during the approximation~\cite{Phan2012,Zhou2012}.

One can impose additional structure on the coefficients of the CP decomposition, such as nonnegativity. As this is of primary interest in data analysis applications, a comprehensive discussion is beyond the scope of this survey, see~\cite{Kolda2009}. 

\subsection{Tucker decomposition} \label{sec:tucker} 

A \emph{Tucker decomposition} of a tensor $\calX$ takes the form
\begin{equation} \label{eq:def_tucker}
\vect(\calX) = 
(U_d \otimes U_{d-1} \otimes \cdots \otimes U_1) \vect(\calC),
\end{equation}
where $U_1, U_2, \ldots,
U_d$, with $U_\mu \in \R^{n_\mu \times r_\mu}$,
are called the \emph{factor matrices} or \emph{basis matrices}
and $\calC \in \R^{r_1 \times \cdots \times r_d}$ is called
the \emph{core tensor} of the decomposition.

Like CP, the Tucker
decomposition has a long history and we refer to the survey~\cite{Kolda2009} for a more detailed account. In the following, we briefly summarize the basic techniques, which are
needed to motivate the TT and HT decompositions discussed below.

The Tucker decomposition is closely related to the matricizations of $\calX$.
The $\mu$th matricization $X^{(\mu)}$ is an $n_\mu \times \big(n_1\cdots n_{\mu-1} n_{\mu+1}\cdots n_d\big)$ matrix
formed in a specific way from the entries of $\calX$:
\[X^{(\mu)}_{i_\mu,\ell} := \calX_{i_1,\ldots,i_d},\qquad
\ell = 1+\sum_{\nu<\mu}(i_\nu-1)\prod_{\eta<\nu} n_\eta + \sum_{\nu>\mu}(i_\nu-1)\prod_{\eta<\nu\atop \eta\ne\mu} n_\eta.\]
In particular, the relation~\eqref{eq:def_tucker} implies
\begin{equation} \label{eq:matrixtucker}
 X^{(\mu)} = U_\mu \cdot C^{(\mu)} \cdot \big( U_d \otimes \cdots \otimes U_{\mu+1} \otimes U_{\mu-1} \otimes \cdots \otimes U_1 \big)^T, \qquad \mu = 1,\ldots, d.
\end{equation}
It follows that $\text{rank} \big( X^{(\mu)} \big) \le r_\mu$, as the first factor
$U_\mu \in \R^{n_\mu\times r_\mu}$ obviously has rank at most $r_\mu$.
This motivates to define the \emph{multilinear rank} (also called $\mu$-rank)
of a tensor $\calX$ as the tuple
\[
(r_1,\ldots, r_d), \quad \text{with}\quad  r_\mu = \text{rank} \big( X^{(\mu)} \big). 
\]
In contrast to the tensor rank related to the CP decomposition, the 
set T($r_1,\ldots,r_d$) of tensors of $\mu$-rank at most $r_\mu$ is closed.

Another consequence of the relation~\eqref{eq:matrixtucker} is the higher-order SVD (HOSVD)
introduced in~\cite{DeLathauwer1997,DeLathauwer2000} for approximating a tensor by a Tucker decomposition~\eqref{eq:def_tucker} of lower multilinear rank. In HOSVD, the columns of each factor
matrix $U_\mu$ are computed as the $k_\mu$ dominant left singular vectors of $X^{(\mu)}$.
The core tensor is then obtained by forming
$
 \text{vec}(\calC):= \big(U_d \otimes \cdots \otimes U_1 \big)^T \text{vec}(\calX).
$
Eventually, this yields 
\[
 {\text{vec}\big( \widetilde \calX\big):= \big(U_d \otimes \cdots \otimes U_1 \big)\cdot \text{vec}(\calC)}
 \quad
 \in T(k_1,\ldots,k_d).
\]
In contrast to the matrix case, where the SVD yields a best low-rank
approximation for all unitarily invariant norms~\cite[Sec. 7.4.9]{Horn2013}, the 
truncated tensor $\widetilde \calX$ resulting from the HOSVD
is usually not optimal.
However, we have
\[
 \textcolor{black}{\| \calX - \widetilde \calX \| \le \sqrt{d} \min_{{\cal Y}\in T(k_1,\ldots,k_d)} \| \calX - {\cal Y}\|}.
 \]
This quasi-optimality condition is usually sufficient for the purpose of obtaining an accurate approximation to a function-related tensor.

Various alternatives to improve on the approximation provided by the HOSVD have been developed,
see~\cite{Kolda2009} and the references therein.
Recent developments include
Newton-type methods on manifolds~\cite{Elden2009,Ishteva2009,Ishteva2011,Savas2010},
a Jacobi algorithm for symmetric tensors~\cite{Ishteva2011a},
generalizations of Krylov subspace methods~\cite{Goreinov2012,Savas2013},
and modifications of the HOSVD~\cite{Vannieuwenhoven2012}.

\subsection{Tensor train decomposition}

The need for storing the $r_1\times \cdots \times r_d$ core tensor $\calC$ renders the Tucker decomposition
increasingly unattractive as $d$ gets larger. This has motivated the search for decompositions which potentially avoid these exponentially growing memory requirements, while still featuring the two most important advantages of the Tucker decomposition: closedness and SVD-based compression.

One well established candidate for such a decomposition is the so called \emph{TT} (tensor train) decomposition,
which takes the form
\begin{equation} \label{eq:tt}
\calX_{i_1,\ldots,i_d} = 
G_1(i_1)\cdot 
G_2(i_2)\cdots
G_d(i_d),\qquad G_\mu(i_\mu) \in\mathbb{R}^{r_{\mu-1} \times r_{\mu}},
\end{equation}
where $r_0=r_d=1$. For every mode $\mu$ and every index $i_{\mu}$ the coefficients $G_{\mu}(i_{\mu})$ are matrices.
In the context of numerical analysis, a decomposition of the form~\eqref{eq:tt} was first
proposed in~\cite{Oseledets2009a,Oseledets2009f,Oseledets2011}.
However, such a decomposition has been proposed earlier in the density-matrix renormalization group method (DMRG) for simulating quantum systems~\cite{White1992,Schollwock2011}. In this area, the term \emph{matrix product state} (MPS) representation for the decomposition~\eqref{eq:tt} has been established~\cite{Ostlund1995}.
Suitable conditions that imply a unique MPS representation can be found in~\cite{Perez-Garcia2007,Vidal2003}.
The connection between TT and MPS has been explained in~\cite{Huckle2012}.

Similar to the Tucker decomposition, the TT decomposition is closely related to
certain matricizations of $\calX$. Let $X^{(1,\ldots,\mu)}$ denote the
matrix obtained by reshaping the entries of $\calX$ into an $(n_1\cdots n_\mu)\times (n_{\mu+1}\cdots n_d)$ array, such that~\eqref{eq:tt} implies $\text{rank}\big(X^{(1,\ldots,\mu)}\big)\le r_{\mu}$ for $\mu = 1,\ldots, d$. Consequently, the tuple containing the ranks of these matricizations is called the \emph{TT-rank} of $\calX$. 
As explained, e.g., in~\cite{Oseledets2011} a quasi-best approximation in a TT decomposition for a given TT-rank can be obtained from the SVDs of $X^{(1,\ldots,\mu)}$, similarly to the HOSVD. It is important to avoid the explicit construction of these matrices and the SVDs when truncating a tensor in TT decomposition to lower TT-rank. Such truncation algorithms are described in~\cite{Oseledets2011}.
On the theoretical side, it turns out that the set $TT(r_1,\ldots,r_{d-1})$ of tensors with TT-ranks bounded 
by $r_\mu$ is closed, and under a full rank condition it 
actually forms a smooth manifold~\cite{Holtz2010,Uschmajew2012}.
The K\"ahler manifold structure for complex MPS with open and periodic boundary
conditions has been studied in~\cite{Haegeman2012}.

\begin{figure}
\begin{center}
  \includegraphics[width=10cm]{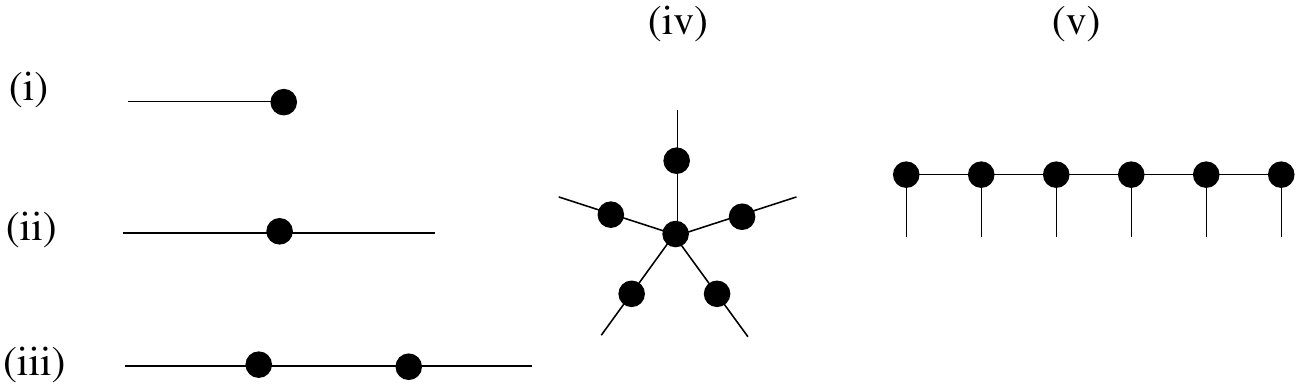}
\end{center}
 \caption{Tensor network diagrams representing
(i) a vector, (ii) a matrix, (iii) a matrix-matrix multiplication, (iv) a tensor in Tucker decomposition,
and (v) a tensor in TT decomposition.} \label{fig:tnexamples}
\end{figure}

Tensor network diagrams, which have been attributed to Penrose~\cite{Penrose1971},
are helpful in visualizing tensor decompositions and their manipulation.
Figure~\ref{fig:tnexamples} gives a few basic examples,
see, e.g.,~\cite{Holtz2012,Huckle2012,Kressner2012} for more details.
In particular, Figure~\ref{fig:tnexamples} (v) gives an illustration of the contraction~\eqref{eq:tt}
representing a TT decomposition. In view of this diagram, the TT decomposition is also sometimes called linear tensor network~\cite{Frowis2010}.

In applications related to quantum spin systems, the tensor $\calX$ often exhibits symmetries inherited from underlying physical properties.
There are variants of MPS/TT that reflect such symmetries in the low-rank decomposition, see~\cite{Huckle2013,Perez-Garcia2007,Sanz2009} and the references therein.

\subsection{Hierarchical Tucker decomposition}

\begin{figure}
\begin{center}
\resizebox{!}{4cm}{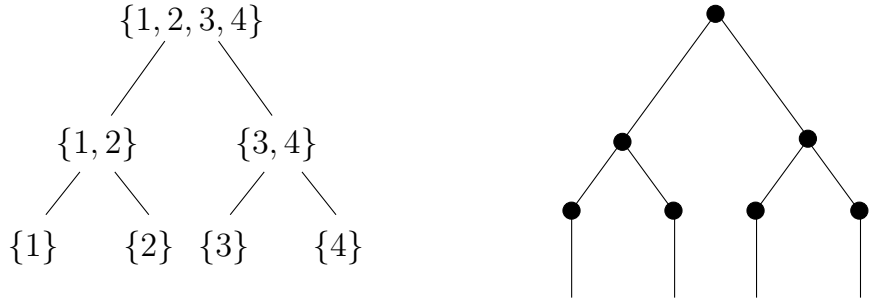}
\end{center}
 \caption{Left: Binary tree representing mode splitting for HT
   decomposition. Right: Tensor network diagram representing a
   tensor in HT decomposition.} \label{fig:tn_ht}
\end{figure}

An alternative way to reduce the complexity of the Tucker decomposition is given by  the \emph{hierarchical Tucker (HT) decomposition} \cite{Hackbusch2009b,Grasedyck2010} (also called
\emph{hierarchical tensor representation}). This decomposition is based on the idea of
recursively splitting the modes of the tensor, which results in a binary
tree $\calT$ containing a subset $t \subset \{1, \ldots, d\}$ at each node. An
example of such a binary tree is given in the left
plot of Figure~\ref{fig:tn_ht}. The matricization $X^{(t)}$ of a tensor $\calX$ corresponding to such a subset $t$ merges all modes contained in $t$ into row indices of the matrix, and the other modes into column indices. We then consider a
hierarchy of matrices $U_t$ whose columns span the image
of $X^{(t)}$ for each $t \in \calT$. Hence, $U_t$ has exactly $r_t =
\text{rank}(X^{(t)})$ columns. The rank tuple $(r_t)_{t\in \calT}$ is called 
the \emph{HT-rank} of $\calX$.

The following nestedness property
allows for the implicit storage of $(U_t)_{t \in \calT}$, and thus of the
tensor $\calX$: For $t = t_l \cup t_r$, $t_l \cap t_r = \emptyset$,
there exists a matrix $B_t$ such that
\begin{equation} \label{eq:htcontained}
U_t = (U_{t_r} \otimes U_{t_l}) B_t,
\qquad
B_t \in \R^{r_{t_l} r_{t_r} \times r_t}.
\end{equation}
For simplicity, we have assumed
that the ordering of the modes in the tree $\calT$ is such that all modes contained in $t_l$ are smaller than the modes contained in $t_r$. The relation~\eqref{eq:htcontained} implies that it suffices to store the basis matrices $U_t$ only for the leaf nodes $t
= \{1\}, \{2\}, \ldots, \{d\}$, and $B_t$ for all other nodes in
$\calT$. The resulting storage requirements are $O(dnr + dr^3)$, when assuming
$r \equiv r_t$ and $n \equiv n_\mu$.

Similarly to the Tucker and TT decompositions, a quasi-best approximation in the HT  decomposition for a given HT-rank can be obtained from the SVDs of $X^{(t)}$. Algorithms that avoid the explicit computation of these SVDs when truncating a tensor that is already in HT decomposition are discussed in~\cite{Grasedyck2010,Hackbusch2009b,Kressner2012,Kuhn2012}.
As for the TT decomposition, the set of tensors having fixed HT-rank forms a smooth manifold~\cite{Falco2013,Uschmajew2012,Uschmajew2013}.

The tensor network corresponding to the HT decomposition is always a binary tree, see also the right plot of Figure~\ref{fig:tn_ht}.
Such tensor tree networks had already been discussed in~\cite{Shi2006}
(without the basis matrices at the leafs). Moreover, the so called
multilayer multi-configuration time-dependent Hartree method
(ML-MCTDH) introduced in~\cite{Wang2003} makes use of a decomposition
based on general trees instead of binary trees. 
When allowing for general trees, tensor tree networks
include the Tucker decomposition from Section~\ref{sec:tucker} as a (quite particular) special case.
In the case of a degenerate tree, where at each level, one mode is
split from the remaining modes, the HT decomposition becomes
equivalent to a \emph{variant} of the TT decomposition discussed
in~\cite{Dolgov2012c,Oseledets2011}. In contrast to the TT
decomposition defined in~\eqref{eq:tt}, this variant features additional basis
matrices, which may reduce the storage cost for large $n_\mu$.
A discussion on the difference between the ranks for the HT and TT decompositions can be
found in~\cite{Grasedyck2011}.

\subsection{More general tensor network formats}

Motivated by an underlying topology describing interactions, tensor networks beyond trees have been
considered in the context of renormalization
group methods for simulating strongly correlated quantum
spin systems. Well-known examples include the so called projected entangled-pair
states (PEPS)~\cite{Verstraete2004a,Verstraete2004b} and
the multiscale entanglement renormalization ansatz (MERA)~\cite{Vidal2007}.
Both, PEPS and MERA contain cycles in the tensor network.
Tree-structured tensor networks, as the hierarchical Tucker and the TT format, are
closed~\cite{Espig2011,Falco2012,Hackbusch2012} in the sense that tensors with ranks at most $r_\mu$ 
form a closed set in $\R^{n_1\times \cdots \times n_d}$.
In general, this statement does not hold for tensor networks containing cycles~\cite{Landsberg2012a,Landsberg2012}.
Possibly for this reason, more general networks have not yet been considered to a large
extend in the numerical analysis community for, e.g., the solution of high-dimensional PDEs,
but see~\cite{Espig2011,Handschuh2012} for some recent mathematically oriented work.

\subsection{Hybrid formats}

Adding to the diversity of the formats discussed above, it is possible and sometimes useful
to combine different low-rank formats. One popular combination is the Tucker format combined with
the CP format for the approximation of the core tensor~\cite{Khoromskij2006,Khoromskij2007a,Khoromskij2009,Kressner2010},
see~\cite[Sec. 5.7]{Kolda2009} for other variations of Tucker and CP.
In~\cite{Gavrilyuk2005,Hackbusch2005,Hackbusch2006,Hackbusch2006a}, combinations of low-rank tensor formats with hierarchical matrices
are investigated.


\subsection{A priori approximation results}

As an essential prerequisite for the success of tensor-based computations, it is important to decide whether a tensor generated by a certain multivariate function $f(x_1,\ldots,x_d)$ can be
well approximated by a low-rank tensor decomposition. As discussed in Section~\ref{sec:cp},
the tensor rank is closely linked to approximating $f$ by a sum of separable functions.
Only in exceptional cases, it will be possible to represent $f$ \emph{exactly} by such a sum,
see~\cite{Beylkin2002,Mohlenkamp2005,Oseledets2012a}.

In general, one is therefore interested in an approximation of the form 
\begin{equation} \label{eq:sumseparable}
 f(x_1,\ldots,x_d) = \sum_{r = 1}^R f^{(1)}_r(x_1)\cdot  f^{(2)}_r(x_2)\cdots f^{(d)}_r(x_d) + \varepsilon.
\end{equation}
For a function of the form $f(x_1,\ldots,x_d) = g(x_1+\cdots+x_d)$,
such an approximation can be immediately obtained from approximating $g$ by a sum of exponentials. For this purpose, various approaches have been discussed in~\cite{Beylkin2005,Beylkin2005a,Beylkin2010,Braess2005,Braess2009}.
Other techniques include applying numerical quadrature to an integral representation
of the function, see, e.g.,~\cite{Bertoglio2012,Espig2013a,Flad2010,Hackbusch2005,Hackbusch2006}, Taylor series expansion~\cite{Tyrtyshnikov2003,Tyrtyshnikov2004},
and polynomial interpolation~\cite{Bebendorf2000,Borm2006}.
Based on results by Temlyakov~\cite{Temlyakov1986,Temlyakov1992,Temlyakov1992a}
on bilinear approximation rates, singular value estimates for the hierarchical Tucker decomposition
of functions in mixed Sobolev spaces have been obtained in~\cite{Schneider2013}, see also~\cite{Griebel2011}.
In fact, a sparse grid approximation to $f$ can be turned quite effectively into 
a low-rank tensor decomposition~\cite{Hackbusch2009b,Hackbusch2012}.
General nonlinear best $R$-term approximation schemes~\cite{Cohen2010,Cohen2011,Kressner2011} represent another important technique, which we cannot cover in detail.


Even for smooth $f$, it may not always be possible to attain sufficiently low ranks, especially when the variation of $f$ is too strong across its entire domain of definition. In this case, it can be advantageous to subdivide the domain and approximate $f$ on each subdomain separately with a low-rank tensor decomposition, see~\cite{Bachmayr2012,Ballani2012c,Ballani2012b} for examples.
As first discussed in~\cite{Beylkin2002}, an approximation of the form~\eqref{eq:sumseparable} can also be used to approximate linear operators on tensors, see also Section~\ref{sec:linop} below.

In the other direction, SVD-based approximations of a function-related tensor yield an approximation of the underlying function, where the $L^2$-norm approximation error can be directly controlled by the truncated singular values.
For smooth functions, best approximations in tensor formats are known to
inherit the regularity of the approximated function~\cite{Uschmajew2011,Uschmajew2013}. This also
holds for SVD based quasi-best approximations, for which even the  smoothness of the error can be
controlled~\cite{Hackbusch2012b}. This can be used, e.g., for deriving $L^\infty$ error estimates~\cite{Hackbusch2012b}
or approximation results for the basis matrices $U_\mu$~\cite{Schneider2013}.

For quantum many-body systems, the approximability of the ground state by a low-rank TT decomposition
is closely linked to the concepts of entropy and entanglement, see~\cite{Perez-Garcia2007} for an introduction.

\subsection{Low-rank decomposition of linear operators} \label{sec:linop}

The matrix representation of a linear operator
\[
 \calA: \R^{n_1\times n_2 \times \cdots \times n_d} \to 
 \R^{m_1\times m_2 \times \cdots \times m_d}, \quad \calX \mapsto \calA(\calX),
\]
can be viewed as an $m_1n_1 \times m_2n_2 \times \cdots \times m_dn_d$
tensor after pairing up row and column indices: 
\[
 \calA_{i_1,i_2,\ldots,i_d;
 j_1,j_2,\ldots,j_d} \qquad \Rightarrow \qquad
 \calA_{(i_1,j_1),(i_2,j_2),\ldots,(i_d,j_d)}.
\]
This view allows to apply any of the low-rank tensor decompositions discussed
above to $\calA$. Such a low-rank decomposition of $\calA$ is useful, e.g., 
for performing the matrix-vector product $\calA(\calX)$
efficiently when $\calX$ itself admits a low-rank decomposition.
This idea appears to be ubiquitous in the literature on low-rank tensor
decomposition, see~\cite{Beylkin2002} for an early reference. For example, in the study of strongly correlated quantum systems,
matrix product operators (MPO) were introduced in~\cite{Verstraete2004,Zwolak2004}, which corresponds to representing $\calA$ in the TT decomposition.
As pointed out in~\cite{Hackbusch2008}, having $\calA$ in a low-rank decomposition also allows to compute an approximate inverse 
by combining the Newton-Hotelling-Schulz algorithm with truncation, see  Section~\ref{sec:iteratetruncate}.

\subsection{Tensorization}

Reverting the process of vectorization,
the entries $x_j$ of a vector $x \in \R^{N}$, $N=2^d$, can be
rearranged into an $2\times 2 \times \cdots \times 2$ tensor
$\calX$ of order $d$. 
(The binary representation $j-1=\sum_{\mu=1}^{d}2^{\mu-1}i_\mu$ yields a simple mapping to the corresponding 
multi-index $(i_1,\ldots,i_d)$ of $\calX$.)
In turn, a low-rank approximation of $\calX$ yields a compression of the
original vector $x$. In combination with the TT decomposition,
this idea of \emph{tensorization} or \emph{quantization}
is usually called \emph{Quantics-TT} or \emph{Quantized-TT} (QTT);
it was first used as a compression scheme for matrices in~\cite{Oseledets2010b}, and introduced for a more general setting in~\cite{Khoromskij2011c}.

Quantization is particularly interesting when the vector $x$ represents a function $f: I \to \R$ evaluated at $2^d$ points, usually uniformly distributed in the interval. The exact and approximate ranks of $\calX$
for various functions $f$ have been discussed for the TT and HT decompositions
in~\cite{Grasedyck2010b,Hackbusch2012,Khoromskij2011c}.

Applying quantization to each mode of a tensor
$\calX \in \R^{N \times \cdots \times N}$ of order $D$, $N=2^d$,
yields a $2\times 2\times \cdots \times 2$ tensor $\calY$ of order $d\cdot D$.
This gives rise to a variety of mixed low-rank tensor decompositions, as discussed in~\cite{Dolgov2012c,Khoromskij2011c}. 

QTT has been applied to the solution of PDEs and eigenvalue problems~\cite{Khoromskij2011a,Khoromskij2011c,Khoromskaia2011}, evaluation of boundary integrals in BEM~\cite{Khoromskij2011f}, convolution~\cite{Hackbusch2011} and the FFT~\cite{Dolgov2012b,Hackbusch2012,Savostyanov2012a}.
A connection between QTT and the wavelet transform is discussed in~\cite{Oseledets2011b}.

One important ingredient of QTT is that the involved
matrices can be represented in a way that conforms to the format~\cite{Khoromskij2011a}. For the following matrices, QTT representations have been discussed: Toeplitz matrices~\cite{Kazeev2011,Oseledets2011d}, (inverse) Laplace operators~\cite{Kazeev2012,Oseledets2010b}, linear diffusion operators~\cite{Dolgov2012c,Kazeev2012a}.

\subsection{Software}

There are several {\sc Matlab} toolboxes available for dealing with tensors in CP and Tucker decomposition, including the Tensor Toolbox~\cite{Bader2010,Bader2006},
the $N$-way toolbox \cite{Andersson2000}, the PLS\_Toolbox~\cite{Wise2013}, and the Tensorlab~\cite{Sorber2013}.
The TT-Toolbox \cite{Oseledets2011a} provides {\sc Matlab} classes
covering tensors in TT and QTT decomposition, as well as linear
operators. There is also a Python implementation of the TT-toolbox called {\tt ttpy}~\cite{Oseledets2013}. The \texttt{htucker} toolbox
\cite{Kressner2012} provides a {\sc Matlab} class representing a tensor in
HT decomposition.

The \emph{TensorCalculus library}~\cite{Espig2012b} is a mathematically oriented C++ library, allowing for computations with general tensor networks.
The Heidelberg MCTDH Package~\cite{Worth2003} is a set of Fortran programs for multi-dimensional quantum dynamics.
ALPS~\cite{Baueretal2011} provides 
C++ libraries for simulating strongly correlated quantum mechanical systems, including DMRG. Block~\cite{Sharma2013} is a C++ implementation of the DMRG algorithms discussed in~\cite{Sharma2012}.
The tensor contraction engine~\cite{Aueretal2006} automatically generates near-optimal code for tensor contractions in many-body electronic structure methods.

\section{Algorithms} \label{sec:algorithms}

In applications for function-related tensors, $\calX$ is often given implicitly, e.g., as the solution to a linear system or eigenvalue problem. There are mainly two different types of approaches to obtain an approximation to $\calX$. A first class of methods is based on combining classical iterative algorithm with repeated low-rank truncation. A second class is based on reformulating the problem at hand as an optimization problem, constraining the admissible set to low-rank tensors, and applying various optimization techniques.

\subsection{Iterative methods combined with truncation}
\label{sec:iteratetruncate}

In principle, any vector iteration for solving a linear algebra problem involving a tensor can be combined with truncation in any of the low-rank decompositions discussed above.
To illustrate the basic principle, let us consider the (preconditioned) Richardson iteration for solving a linear system $\calA(\calX) = \calB$:
\[
 \calX_{k+1} = \calX_k + \omega \calP\big( \calB - \calA(\calX_k) \big),
\]
where $\calP$ is a preconditioner and $\omega$ is a suitably chosen scalar.
Letting $\calT$ denote truncation in any of the low-rank tensor decompositions
discussed above, we obtain the truncated 
Richardson iteration
\[
 \calX_{k+1} = \calT\Big(\calX_k + \omega \calP\big( \calB - \calA(\calX_k) \big) \Big),
\]
which has been proposed~\cite{Khoromskij2011b} in combination with the 
CP decomposition.

Other examples of combining iterative methods with low-rank truncation include:
\begin{itemize}
 \item 
 the (shift-and-invert) power method combined with CP~\cite{Beylkin2002,Beylkin2005a};
 \item the (shift-and-invert) power and Lanczos methods combined with CP and Tucker~\cite{Bachmayr2012,Hackbusch2012a};
 \item a restarted Lanczos method combined with TT~\cite{Huckle2012a};
 \item conjugate gradient type-methods for symmetric eigenvalue problems combined with truncation
for low-rank matrices~\cite{Lebedeva2010}, as well as for TT~\cite{Mach2011}, HT~\cite{Kressner2011a} and QTT~\cite{Lebedeva2011}.
 \item the Richardson method combined with QTT~\cite{Khoromskij2010} and low-rank matrix decompositions~\cite{Matthies2012};
 \item a projection method combined with HT~\cite{Ballani2012};
 \item the conjugate gradient method, BiCGStab and other Krylov subspace methods combined with HT~\cite{Kressner2011,Kressner2012a,Tobler2012} and low-rank matrix decompositions~\cite{Bollhofer2012};
 \item GMRES combined with TT~\cite{Dolgov2012}.
\end{itemize}
When applying an iterative method in combination with a low-rank decompositions, the ranks will inevitably grow quite quickly in the course of the iteration.
For example, the sum of $k$ tensors can multiply the ranks by $k$, while the pointwise
(Hadamard) product between two tensors can even square the ranks. To gain efficiency, it is therefore advisable to not let this rank growth happen explicitly and combine such an operation directly with truncation. This has been discussed for sums in~\cite{Kressner2012,Tobler2012} and for the Hadamard product (and other bilinear operations) in~\cite{Savostyanov2012}, see also~\cite{Oseledets2011c}.

Preconditioners not only help accelerate convergence but numerical evidence suggests that an effective preconditioner leads to a limited intermediate rank growth.
At the same time, it is mandatory that the preconditioner can be applied efficiently to a low-rank tensor decomposition. In particular, this is the case when the preconditioner itself admits a low-rank decomposition in the sense of Section~\ref{sec:linop}. There are various techniques to construct such 
preconditioners. An early technique is based on best approximation in the Frobenius norm by a Kronecker product~\cite{Ullmann2010,VanLoan1993}, by a short sum of Kronecker products~\cite{Olshevsky2006}, or
by a more general low-rank tensor decomposition. Other techniques include the use of approximate inverses for high-dimensional Laplace operators~\cite{Khoromskij2009b,Khoromskij2011b,Kressner2011a,Oseledets2009b},
low-rank manipulation of the PDE coefficients~\cite{Dolgov2012d},
low-rank tensor approximation of multilevel preconditioners~\cite{Andreev2012b},
and low-rank tensor diagonal preconditioners for wavelet discretizations~\cite{Bachmayr2012}.

%
%

\subsection{Optimization-based algorithms} \label{sec:optimization}

In many cases, a linear algebra problems involving a tensor can be posed as an optimization problem. For example, it is well known that a symmetric positive definite linear system $\calA(\calX) = \calB$ can be turned into
\begin{equation} \label{eq:optimization}
 \min_{\calX \in \R^{n_1\times \cdots \times n_d}} \frac12 \langle \calX, \calA(\calX) \rangle - 
 \langle \calX, \calB \rangle,
\end{equation}
where $\langle \cdot,\cdot \rangle$ corresponds to the standard inner product for the vectorization of the tensors.
For nonsymmetric linear systems, an optimization problem can be obtained by minimizing the norm of the residual. For symmetric eigenvalue problems, the Rayleigh-quotient minimization or, more generally, the trace minimization principle can be used. 

Once the optimization problem is set up,
the set of admissible tensors $\calX$ is then constrained to a low-rank decomposition, for example, to all tensors with fixed tensor rank or fixed multilinear rank. Even when the original optimization problem, such as~\eqref{eq:optimization}, is convex and in principle simple, the resulting constrained optimization 
problem is highly nonlinear and non-convex in general.
A number of heuristic approaches to the solution of such constrained optimization 
problems are available, including ALS (alternating linear scheme). The basic principle of ALS is to optimize every factor of the low-rank decomposition separately and to sweep over all factors repeatedly. It is probably most natural to combine ALS with the CP decomposition, see~\cite{Beylkin2005a,Doostan2009} for an application of this idea to linear systems.
The combination of ALS with the TT decomposition has been considered in~\cite{Dolgov2012e,Dolgov2013,Dolgov2013a,Holtz2012}. The convergence of ALS for the TT decomposition has been studied in~\cite{Rohwedder2011}.

The ALS scheme can be improved in various ways. One quite successful improvement for decompositions described by tensor networks is to join two neighboring factors, optimize the resulting supernode, and split the result
into separate factors by a low-rank factorization. 
Originally, this so called DMRG method had been developed for the simulation of strongly correlated quantum lattice systems, see~\cite{Schollwock2011} for an overview. Later on, the ideas of DMRG have been picked up and extended to other applications in the numerical analysis community in a series of papers~\cite{Dolgov2012e,Holtz2012,Huckle2012,Khoromskij2010a,Kressner2011,Oseledets2011c}.

There is a growing interest in applying so called Riemannian optimization techniques~\cite{Absil2008} to~\eqref{eq:optimization}. Examples
include nonlinear conjugate gradient or Newton-like methods on manifolds of low-rank matrices~\cite{Boumal2011,Meyer2011,Ngo2012,Vandereycken2010,Vandereycken2012} or low-rank tensors~\cite{Curtef2007,Curtef2012,Ishteva2009,Ishteva2011,Uschmajew2012}, see also Section~\ref{sec:dynamical}. For tensors in CP decomposition, the approximate solution of~\eqref{eq:optimization} by gradient techniques has been discussed in~\cite{Espig2012a}.

\subsection{Successive rank-1 approximation}

A tempting and surprisingly successful approach to the solution of high-dimensional problems is
to build up a low-rank approximation from 
successive rank-1 approximations. This idea has been suggested in the context
of various applications, including
Fokker-Planck equations~\cite{Ammar2006} 
and stochastic partial differential equations~\cite{Nouy2010}.

We will illustrate the basic idea with a simple example. Consider
a linear system
$
 \calA(\calX) = \calB
$
with the solution tensor $\calX \in \R^{n_1\times \cdots \times n_d}$. Assume we already have
a CP approximation $\calX_r$ of tensor rank $r$:
\begin{equation}  \label{eq:approx}
 \vect(\calX_r) = 
u^{(d)}_1 \otimes u^{(d-1)}_1 \otimes \cdots \otimes u^{(1)}_1
+ \cdots +
u^{(d)}_r \otimes u^{(d-1)}_r \otimes \cdots \otimes u^{(1)}_r.
\end{equation}
We then search for a rank-1 correction
\[
 \calW = u^{(d)}_{r+1} \otimes u^{(d-1)}_{r+1} \otimes \cdots \otimes u^{(1)}_{r+1},
\]
such that $\calX_{r+1} = \calX_r + \calW$ is an improved approximation, that is,
\begin{equation} \label{eq:correctionequation}
\calA(\calX_{r+1}) \approx \calB \qquad \Leftrightarrow \qquad \calA(\calW) \approx \calB - \calA(\calX_r)
\end{equation}
Analogous to Section~\ref{sec:optimization}, the unknown vectors $w^{(1)}, \ldots, w^{(d)}$ can be determined by turning~\eqref{eq:correctionequation} into a nonlinear (optimization) problem and
applying standard methods, such as the alternating direction method. 
This procedure is repeated until the residual $\calA(\calX_r) - \calB$ is sufficiently small.
Of course, such a greedy approach will not yield the best rank-$R$ approximation after $R$ steps~\cite{Stegeman2010}. However, it is important to remember that these methods aim at a more moderate goal,
to obtain a reasonable approximation after $R$ steps, with $R$ not too large.
Convergence results in this direction can be found in~\cite{Ammar2010,Cances2011,Cances2012,Falco2011,Falco2012a,Figueroa2012,LeBris2009}.

A number of improvements to the simple scheme outlined above have been proposed to increase
its convergence speed, see, e.g.,~\cite{Ammar2006,Giraldi2012,Nouy2007,Nouy2010}.
A connection between best rank-$1$ or, more generally, best rank-$m$ approximations to a
nonlinear eigenvalue problem is explained in~\cite{Nouy2008}. This connection also motivates the
use of the term \emph{generalized spectral decomposition}.
Many further developments, improvements, and extensions of 
successive low-rank approximation techniques have taken place during the last years; we refer to~\cite{Chinesta2011} for an overview.

\subsection{Low-rank methods for dynamical problems} \label{sec:dynamical}

Let us consider a dynamical system on $\R^{n_1\times \cdots \times n_d}$:
\begin{equation} \label{eq:tensorode}
\dot \calX(t) = F(\calX(t)), \quad \calX(t_0) = \calX_0,
\end{equation}
for which a typical application is the spatial discretization of a time-dependent $d$-dimensional PDE.
\emph{Dynamical low-rank methods} aim to determine an approximation $\calY(t)$ in a manifold $\calM$ of low-rank tensors by restricting the dynamics of~\eqref{eq:tensorode} to the tangent space $T_{\calY(t)}\calM$:
\begin{equation} \label{eq:tensorodetangent}
 \dot \calY(t) \in T_{\calY(t)}\calM \quad \text{such that} \quad \|\calY(t) - F(\calY(t))\| = \min!
\end{equation}
As explained in~\cite{Lubich2005}, this approximation is closely related to the  can Dirac-Frenkel-McLachlan variational principle in quantum molecular dynamics.

Initially proposed for low-rank matrix manifolds in~\cite{Koch2007}, dynamical low-rank methods
have been extended to low-rank tensors in Tucker~\cite{Koch2010,Nonnenmacher2008}, TT/MPS~\cite{Haegeman2011,Holtz2010,Khoromskij2012a,Lubich2012}, and HT~\cite{Arnold2012,Lubich2012,Uschmajew2012} decomposition.
The efficient and robust numerical integration of~\eqref{eq:tensorodetangent}
is crucial to the success of dynamical low-rank methods; apart from the references
above, this aspect has been discussed in~\cite{Lubich2008,Lubich2013}.

A more immediate approach to~\eqref{eq:tensorode} is to combine a standard time stepping method, such as the explicit and implicit Euler methods, with low-rank truncation in every time step~\cite{Dolgov2012f}. An alternative, which allows to control the error global-in-time, is to apply iterative solvers to a space-time formulation~\cite{Ammar2007,Andreev2012b,Dolgov2012f,Dolgov2012g,Gavrilyuk2011}.

\subsection{Black box approximation}

Suppose that a matrix $A$ or a tensor $\calX$ is defined through a function that returns 
entries at arbitrary positions. Then the goal of \emph{black box approximation}
is to find a good low-rank approximation based only on relatively few entries.
It is important to emphasize that the selection of the entries can be controlled by the
user. In this respect, this situation is quite different from the growing area of tensor completion, see, e.g.,~\cite{Gandy2011,Liu2009}, where the selection of the entries is usually prescribed by the application.

For an $m\times n$ matrix $A$, the so called \emph{cross approximation} method~\cite{Bebendorf2000,Borm2005,Goreinov1997,Tyrtyshnikov2000} 
produces an approximation of the form
\begin{equation} \label{eq:crossapproximation}
 A(:,J) A(I,J)^{-1} A(I,:),
\end{equation}
where {\sc Matlab} notation is used
to denote the submatrices of $A$ corresponding to the index sets
$I \subset \{1,\ldots,m\}$ and $J \subset \{1,\ldots,n\}$ .
In the $p$th step of cross approximation as described in~\cite{Bebendorf2003,Oseledets2008}, the entry of largest magnitude in the column $j_{p}$
of $A - A(:,J) A(I,J)^{-1} A(I,:)$ is calculated and its position is
denoted by $i_p$. Then the entry of largest magnitude in the row $i_p$
of that matrix is calculated and its position is denoted by 
$j_{p+1}$. Moreover, both index sets are updated: $I\leftarrow I \cup \{i_p\}$
and $J\leftarrow J \cup \{j_p\}$.
Volume maximization is an alternative entry selection strategy that attempts to maximize
$\big|\det\big(A(I,J)\big)\big|$, see~\cite{Goreinov2010,Tyrtyshnikov2000}.

A first extension of cross approximation to tensors was proposed in~\cite{Oseledets2008}
for approximating third-order tensors by a Tucker decomposition. Essentially, this extension consists of applying the algorithm above to an arbitrary matricization of the tensor. However, the rows of this matricization (corresponding to slices of the tensor) are further approximated by a low rank matrix, again using cross approximation.
In~\cite{Flad2008,Oseledets2010c}, this method has been combined with multilevel ideas and applied to quantum chemistry.
The adaptive cross approximation for the Tucker decomposition has been analyzed in~\cite{Bebendorf2011}.
Another variant for third and fourth order tensors, focussing on interpolation properties, is discussed in~\cite{Naraparaju2011}.
A cross approximation method for the TT decomposition has been proposed in~\cite{Oseledets2010a,Savostyanov2011a}.

Based on fiber crosses, a quite different extension for tensors of arbitrary order can be found in~\cite{Espig2009}.
In this method a set of multi-indices $I^1,\ldots,I^p \in [1,n_1]\times \cdots [1,n_d]$ is computed successively. Subspaces $\calU_\mu$ are constructed approximately containing all $\mu$-mode fibers passing through at least one of the multi-indices. The tensor is then approximated
by a CP decomposition~\eqref{eq:def_cp} under the constraint $u_j^{(\mu)} \in \calU_\mu$, using general
optimization methods. The multi-indices are selected by performing an alternating
direction search along the fibers of the error tensor.
Also based on fiber crosses, a black box approximation for a tensor in the HT decomposition is given in~\cite{Ballani2012a,Ballani2012c}.

Randomized algorithms represent an alternative way to quickly extract a low-rank approximation from partial information on the entries of the  matrix, see~\cite{Halko2011} and the references therein. These ideas have been
extended to low-rank tensor decompositions, for the Tucker decomposition in~\cite{Drineas2007,Friedland2011,Tsourakakis2010} and for the TT decomposition in~\cite{Pan2012}.

\subsection{Other algorithms}

For linear systems and eigenvalue problems with a very particular structure, it is sometime
possible to design specialized algorithms that can be more efficient and easier to analyze.
This applies, in particular, to discretizations of the multi-dimensional Poisson equation~\cite{Beckermann2012,Grasedyck2004,Kressner2010,Mach2012}.

\section{Acknowledgments}

We thank Jonas Ballani, Wolfgang Hackbusch, Thomas Huckle, Boris Khoromskij, Anthony Nouy, Ivan Oseledets, Andr\'e Uschmajew, and
Bart Vandereycken for many helpful suggestions on a preliminary draft of this survey.

\bibliographystyle{gamm}
\bibliography{tensor_refs}

\providecommand{\WileyBibTextsc}{}
\let\textsc\WileyBibTextsc
\providecommand{\othercit}{}
\providecommand{\jr}[1]{#1}
\providecommand{\etal}{~et~al.}


\begin{thebibliography}{[100]}

\othercit
\bibitem{Absil2008}
 \textsc{P.\,A. Absil},  \textsc{R.~Mahony},  and  \textsc{R.~Sepulchre},
Optimization algorithms on matrix manifolds (Princeton University Press,
  Princeton, NJ, 2008).


\bibitem{Acar2011}
 \textsc{E.~Acar},  \textsc{D.\,M. Dunlavy},  and  \textsc{T.\,G. Kolda},
A scalable optimization approach for fitting canonical tensor decompositions,
 \jr{J. Chemometrics} \textbf{25}(2), 67--86 (2011).


\bibitem{Ammar2010}
 \textsc{A.~Ammar},  \textsc{F.~Chinesta},  and  \textsc{A.~Falc{\'o}},
On the convergence of a greedy rank-one update algorithm for a class of linear
  systems,
 \jr{Arch. Comput. Methods Eng.} \textbf{17}(4), 473--486 (2010).


\bibitem{Ammar2006}
 \textsc{A.~Ammar},  \textsc{B.~Mokdad},  \textsc{F.~Chinesta},  and
  \textsc{R.~Keunings},
A new family of solvers for some classes of multidimensional partial
  differential equations encountered in kinetic theory modeling of complex
  fluids,
 \jr{J. Non-Newton. Fluid Mech.} \textbf{139}(3), 153--176 (2006).


\bibitem{Ammar2007}
 \textsc{A.~Ammar},  \textsc{B.~Mokdad},  \textsc{F.~Chinesta},  and
  \textsc{R.~Keunings},
A new family of solvers for some classes of multidimensional partial
  differential equations encountered in kinetic theory modelling of complex
  fluids: Part {II}: Transient simulation using space-time separated
  representations,
 \jr{J. Non-Newton. Fluid Mech.} \textbf{144}(2--3), 98--121 (2007).


\bibitem{Ammar2010a}
 \textsc{A.~Ammar},  \textsc{M.~Normandin},  \textsc{F.~Daim},
  \textsc{D.~Gonz{\'a}lez},  \textsc{E.~Cueto},  and  \textsc{F.~Chinesta},
Non incremental strategies based on separated representations: applications in
  computational rheology,
 \jr{Commun. Math. Sci.} \textbf{8}(3), 671--695 (2010).


\bibitem{Andersson2000}
 \textsc{C.\,A. Andersson} and  \textsc{R.~Bro},
The {$N$-way} toolbox for {MATLAB},
 \jr{Chemometrics and Intelligent Laboratory Systems} \textbf{52}(1), 1--4
  (2000),
Available at \url{http://www.models.life.ku.dk/nwaytoolbox}.


\othercit
\bibitem{Andreev2012b}
 \textsc{R.~Andreev} and  \textsc{C.~Tobler},
Multilevel preconditioning and low rank tensor iteration for space-time
  simultaneous discretizations of parabolic {PDE}s,
Tech. Rep. 2012--16, Seminar for Applied Mathematics, ETH Z{\"u}rich, 2012,
Available at \url{http://www.sam.math.ethz.ch/reports/2012/16}.


\othercit
\bibitem{Arnold2012}
 \textsc{A.~Arnold} and  \textsc{T.~Jahnke},
On the approximation of high-dimensional differential equations in the
  hierarchical {T}ucker format,
Tech. rep., KIT, Karlsruhe, Germany, 2012.


\bibitem{Aueretal2006}
 \textsc{A.\,A. {Auer et al.}},
Automatic code generation for many-body electronic structure methods: the
  tensor contraction engine,
 \jr{Molecular Physics} \textbf{104}(2), 211--228 (2006).


\othercit
\bibitem{Bachmayr2012}
 \textsc{M.~Bachmayr},
Adaptive low-rank wavelet methods and applications to two-electron
  {S}chr\"odinger equations,
PhD thesis, RWTH Aachen, Germany, 2012.


\bibitem{Bader2006}
 \textsc{B.\,W. Bader} and  \textsc{T.\,G. Kolda},
Algorithm 862: {MATLAB} tensor classes for fast algorithm prototyping,
 \jr{ACM Trans. Math. Software} \textbf{32}(4), 635--653 (2006),
Available at \url{http://csmr.ca.sandia.gov/~tgkolda/TensorToolbox/}.


\othercit
\bibitem{Bader2010}
 \textsc{B.\,W. Bader} and  \textsc{T.\,G. Kolda},
{MATLAB Tensor Toolbox Version 2.4}, March 2010,
Available at \url{http://csmr.ca.sandia.gov/~tgkolda/TensorToolbox/}.


\othercit
\bibitem{Ballani2012c}
 \textsc{J.~Ballani},
{Fast evaluation of near-field boundary integrals using tensor approximations},
Dissertation, Universit\"at Leipzig, 2012.


\bibitem{Ballani2012b}
 \textsc{J.~Ballani},
Fast evaluation of singular {BEM} integrals based on tensor approximations,
 \jr{Numer. Math.} \textbf{121}(3), 433--460 (2012).


\bibitem{Ballani2012}
 \textsc{J.~Ballani} and  \textsc{L.~Grasedyck},
A projection method to solve linear systems in tensor format,
 \jr{Numer. Linear Algebra Appl.} \textbf{20}(1), 27--43 (2013).


\bibitem{Ballani2012a}
 \textsc{J.~Ballani},  \textsc{L.~Grasedyck},  and  \textsc{M.~Kluge},
Black box approximation of tensors in hierarchical {T}ucker format,
 \jr{Linear Algebra Appl.} \textbf{438}(2), 639--657 (2013).


\bibitem{Baueretal2011}
 \textsc{B.~{Bauer et al.}},
The {ALPS} project release 2.0: open source software for strongly correlated
  systems,
 \jr{J. Stat. Mech.} \textbf{5} (2011),
Available at \url{http://alps.comp-phys.org}.


\bibitem{Bebendorf2000}
 \textsc{M.~Bebendorf},
Approximation of boundary element matrices,
 \jr{Numer. Math.} \textbf{86}(4), 565--589 (2000).


\bibitem{Bebendorf2011}
 \textsc{M.~Bebendorf},
Adaptive cross approximation of multivariate functions,
 \jr{Constr. Approx.} \textbf{34}, 149--179 (2011).


\bibitem{Bebendorf2003}
 \textsc{M.~Bebendorf} and  \textsc{S.~Rjasanow},
Adaptive low-rank approximation of collocation matrices,
 \jr{Computing} \textbf{70}(1), 1--24 (2003).


\othercit
\bibitem{Beckermann2012}
 \textsc{B.~Beckermann},  \textsc{D.~Kressner},  and  \textsc{C.~Tobler},
An error analysis of {G}alerkin projection methods for linear systems with
  tensor product structure,
Tech. rep., MATHICSE, EPF Lausanne, Switzerland, 2012.


\bibitem{Benedikt2011}
 \textsc{U.~Benedikt},  \textsc{A.\,A. Auer},  \textsc{M.~Espig},  and
  \textsc{W.~Hackbusch},
Tensor decomposition in post-{H}artree-{F}ock methods. {I}. two-electron
  integrals and {MP2},
 \jr{J. Chem. Phys} \textbf{134}(5), 054118--1--054118--12 (2011).


\bibitem{Bertoglio2012}
 \textsc{C.~Bertoglio} and  \textsc{B.\,N. Khoromskij},
Low-rank quadrature-based tensor approximation of the {G}alerkin projected
  {Newton/Yukawa} kernels,
 \jr{Computer Physics Communications} \textbf{183}(4), 904--912 (2012).


\bibitem{Beylkin2009}
 \textsc{G.~Beylkin},  \textsc{J.~Garcke},  and  \textsc{M.\,J. Mohlenkamp},
Multivariate regression and machine learning with sums of separable functions,
 \jr{SIAM J. Sci. Comput.} \textbf{31}(3), 1840--1857 (2009).


\bibitem{Beylkin2002}
 \textsc{G.~Beylkin} and  \textsc{M.\,J. Mohlenkamp},
Numerical operator calculus in higher dimensions,
 \jr{Proc. Natl. Acad. Sci. USA} \textbf{99}(16), 10246--10251 (2002).


\bibitem{Beylkin2005a}
 \textsc{G.~Beylkin} and  \textsc{M.\,J. Mohlenkamp},
Algorithms for numerical analysis in high dimensions,
 \jr{SIAM J. Sci. Comput.} \textbf{26}(6), 2133--2159 (2005).


\bibitem{Beylkin2005}
 \textsc{G.~Beylkin} and  \textsc{L.~Monz{\'o}n},
On approximation of functions by exponential sums,
 \jr{Appl. Comput. Harmon. Anal.} \textbf{19}(1), 17--48 (2005).


\bibitem{Beylkin2010}
 \textsc{G.~Beylkin} and  \textsc{L.~Monz{\'o}n},
Approximation by exponential sums revisited,
 \jr{Appl. Comput. Harmon. Anal.} \textbf{28}(2), 131--149 (2010).


\othercit
\bibitem{Biagioni2012}
 \textsc{D.\,J. Biagioni},
Numerical construction of {G}reen's functions in high dimensional elliptic
  problems with variable coefficients and analysis of renewable energy data via
  sparse and separable approximations,
PhD thesis, University of Colorado, 2012.


\bibitem{Blesgen2012}
 \textsc{T.~Blesgen},  \textsc{V.~Gavini},  and  \textsc{V.~Khoromskaia},
Approximation of the electron density of aluminium clusters in tensor-product
  format,
 \jr{J. Comput. Phys.} \textbf{231}(6), 2551--2564 (2012).


\othercit
\bibitem{Bollhofer2012}
 \textsc{M.~Bollh\"{o}fer} and  \textsc{A.\,K. Eppler},
Low-rank {C}holesky-factor {K}rylov subspace methods for generalized projected
  {L}yapunov equations,
Tech. rep., To appear in "System Reduction for Nanoscale IC Design" (Peter
  Benner, ed.), Mathematics in Industry (Springer), 2012.


\othercit
\bibitem{Borm2006}
 \textsc{S.~B{\"o}rm},
{$\mathcal{H}_2$}-matrices -- an efficient tool for the treatment of dense
  matrices,
Habilitationsschrift, Christian-Albrechts-Universit\"at zu Kiel, 2006.


\bibitem{Borm2005}
 \textsc{S.~B{\"o}rm} and  \textsc{L.~Grasedyck},
Hybrid cross approximation of integral operators,
 \jr{Numer. Math.} \textbf{101}(2), 221--249 (2005).


\othercit
\bibitem{Boumal2011}
 \textsc{N.~Boumal} and  \textsc{P.\,A. Absil},
{RTRMC}: A {R}iemannian trust-region method for low-rank matrix completion,
 in: Proceedings of the Neural Information Processing Systems Conference
  (NIPS),  (2011).


\bibitem{Braess2005}
 \textsc{D.~Braess} and  \textsc{W.~Hackbusch},
Approximation of {$1/x$} by exponential sums in {$[1,\infty)$},
 \jr{IMA J. Numer. Anal.} \textbf{25}(4), 685--697 (2005).


\othercit
\bibitem{Braess2009}
 \textsc{D.~Braess} and  \textsc{W.~Hackbusch},
On the efficient computation of high-dimensional integrals and the
  approximation by exponential sums,
 in: Multiscale, Nonlinear and Adaptive Approximation,  (Springer, 2009),
  pp.\,39--74.


\bibitem{Buchholz2010}
 \textsc{P.~Buchholz},
Product form approximations for communicating {M}arkov processes,
 \jr{Performance Evaluation} \textbf{67}(9), 797--815 (2010).


\bibitem{Cances2011}
 \textsc{E.~Canc\`es},  \textsc{V.~Ehrlacher},  and  \textsc{T.~Leli\`evre},
Convergence of a greedy algorithm for high-dimensional convex nonlinear
  problems,
 \jr{Mathematical Models and Methods in Applied Sciences} \textbf{21}(12),
  2433--2467 (2011).


\othercit
\bibitem{Cances2012}
 \textsc{E.~Canc\`es},  \textsc{V.~Ehrlacher},  and  \textsc{T.~Leli\`evre},
Greedy algorithms for high-dimensional non-symmetric linear problems,
arxiv:1210.6688, 2012.


\bibitem{Chen2011}
 \textsc{Y.~Chen},  \textsc{D.~Han},  and  \textsc{L.~Qi},
New {ALS} methods with extrapolating search directions and optimal step size
  for complex-valued tensor decompositions,
 \jr{Signal Processing, IEEE Transactions on} \textbf{59}(12), 5888--5898
  (2011).


\bibitem{Chinesta2008}
 \textsc{F.~Chinesta},  \textsc{A.~Ammar},  \textsc{F.~Lemarchand},
  \textsc{P.~Beauchene},  and  \textsc{F.~Boust},
Alleviating mesh constraints: Model reduction, parallel time integration and
  high resolution homogenization,
 \jr{Comput. Meth. Appl. Mech. Engrg.} \textbf{197}(5), 400--413 (2008).


\bibitem{Chinesta2011}
 \textsc{F.~Chinesta},  \textsc{P.~Ladeveze},  and  \textsc{E.~Cueto},
A short review on model order reduction based on proper generalized
  decomposition,
 \jr{Arch. Comput. Methods Eng.} \textbf{18}, 395--404 (2011).


\bibitem{Chinnamsetty2010}
 \textsc{S.\,R. Chinnamsetty},  \textsc{M.~Espig},  \textsc{H.\,J. Flad},  and
  \textsc{W.~Hackbusch},
Canonical tensor products as a generalization of {Gaussian}-type orbitals,
 \jr{Z. Phys. Chem.} \textbf{224}(3--4), 681--694 (2010).


\bibitem{Chinnamsetty2007}
 \textsc{S.\,R. Chinnamsetty},  \textsc{M.~Espig},  \textsc{B.\,N. Khoromskij},
   \textsc{W.~Hackbusch},  and  \textsc{H.\,J. Flad},
Tensor product approximation with optimal rank in quantum chemistry,
 \jr{J. Chem. Phys} \textbf{127}(8), 084110 (2007).


\othercit
\bibitem{Chinnamsetty2009}
 \textsc{S.\,R. Chinnamsetty},  \textsc{W.~Hackbusch},  and  \textsc{H.\,J.
  Flad},
The tensor product approximation to single-electron systems,
Tech. Rep.~23, MPI MIS Leipzig, 2009.


\bibitem{Chinnamsetty2010a}
 \textsc{S.\,R. Chinnamsetty},  \textsc{W.~Hackbusch},  and  \textsc{H.\,J.
  Flad},
Efficient multi-scale computation of products of orbitals in electronic
  structure calculations,
 \jr{Comput. Vis. Sci.} \textbf{13}(8), 397--408 (2010).


\bibitem{Chinnamsetty2012}
 \textsc{S.\,R. Chinnamsetty},  \textsc{H.~Luo},  \textsc{W.~Hackbusch},
  \textsc{H.\,J. Flad},  and  \textsc{A.~Uschmajew},
Bridging the gap between quantum {Monte Carlo} and {F12}-methods,
 \jr{J. Chem. Phys} \textbf{401}, 36--44 (2012).


\bibitem{Cohen2010}
 \textsc{A.~Cohen},  \textsc{R.~DeVore},  and  \textsc{C.~Schwab},
Convergence rates of best {$N$}-term {G}alerkin approximations for a class of
  elliptic s{PDE}s,
 \jr{Found. Comput. Math.} \textbf{10}(6), 615--646 (2010).


\bibitem{Cohen2011}
 \textsc{A.~Cohen},  \textsc{R.~DeVore},  and  \textsc{C.~Schwab},
Analytic regularity and polynomial approximation of parametric and stochastic
  elliptic {PDE}'s,
 \jr{Anal. Appl. (Singap.)} \textbf{9}(1), 11--47 (2011).


\bibitem{Comon2009}
 \textsc{P.~Comon},  \textsc{X.~Luciani},  and  \textsc{A.\,L.\,F. de~Almeida},
Tensor decompositions, alternating least squares and other tales,
 \jr{J. Chemometrics} \textbf{23}(7--8), 393--405 (2009).


\bibitem{Curtef2007}
 \textsc{O.~Curtef},  \textsc{G.~Dirr},  and  \textsc{U.~Helmke},
Conjugate gradient algorithms for best rank-1 approximation of tensors,
 \jr{PAMM} \textbf{7}(1), 1062201--1062202 (2007).


\bibitem{Curtef2012}
 \textsc{O.~Curtef},  \textsc{G.~Dirr},  and  \textsc{U.~Helmke},
{R}iemannian optimization on tensor products of {G}rassmann manifolds:
  Applications to generalized {R}ayleigh-quotients,
 \jr{SIAM J. Matrix Anal. Appl.} \textbf{33}(1), 210--234 (2012).


\othercit
\bibitem{DeLathauwer1997}
 \textsc{L.~De~Lathauwer},
Signal Processing Based on Multilinear Algebra,
PhD thesis, Katholike Universiteit Leuven, 1997.


\bibitem{DeLathauwer2000}
 \textsc{L.~De~Lathauwer},  \textsc{B.~De~Moor},  and  \textsc{J.~Vandewalle},
A multilinear singular value decomposition,
 \jr{SIAM J. Matrix Anal. Appl.} \textbf{21}(4), 1253--1278 (2000).


\bibitem{Silva2008}
 \textsc{V.~de~Silva} and  \textsc{L.\,H. Lim},
Tensor rank and the ill-posedness of the best low-rank approximation problem,
 \jr{SIAM J. Matrix Anal. Appl.} \textbf{20}(3), 1084--1127 (2008).


\bibitem{Sterck2012}
 \textsc{H.~de~Sterck},
A nonlinear {GMRES} optimization algorithm for canonical tensor decomposition,
 \jr{SIAM J. Sci. Comput.} \textbf{34}(3), A1351--A1379 (2012).


\othercit
\bibitem{DeSterck2011}
 \textsc{H.~De~Sterck} and  \textsc{K.~Miller},
An adaptive algebraic multigrid algorithm for low-rank canonical tensor
  decomposition,
Tech. Rep. arXiv:1111.6091, Nov 2011.


\othercit
\bibitem{Dolgov2012}
 \textsc{S.\,V. Dolgov},
{TT-GMRES:} on solution to a linear system in the structured tensor format,
{arXiv} preprint 1206.5512 (To appear in: {Rus. J. of Num. An. and Math.
  Model.}), 2012.


\othercit
\bibitem{Dolgov2012d}
 \textsc{S.\,V. Dolgov},  \textsc{V.\,A. Kazeev},  and  \textsc{B.\,N.
  Khoromskij},
The tensor-structured solution of one-dimensional elliptic differential
  equations with high-dimensional parameters,
Tech. Rep.~51, MPI MIS Leipzig, 2012.


\othercit
\bibitem{Dolgov2012g}
 \textsc{S.\,V. Dolgov} and  \textsc{B.\,N. Khoromskij},
Tensor-product approach to global time-space-parametric discretization of
  chemical master equation,
Tech. Rep.~68, MPI MIS Leipzig, 2012.


\othercit
\bibitem{Dolgov2012c}
 \textsc{S.\,V. Dolgov} and  \textsc{B.\,N. Khoromskij},
Two-level {Tucker-TT-QTT} format for optimized tensor calculus,
Tech. Rep.~19, MPI MIS Leipzig, 2012.


\bibitem{Dolgov2012b}
 \textsc{S.\,V. Dolgov},  \textsc{B.\,N. Khoromskij},  and  \textsc{D.\,V.
  Savostyanov},
Superfast {Fourier} transform using {QTT} approximation,
 \jr{J. Fourier Anal. Appl.} \textbf{18}(5), 519--953 (2012).


\bibitem{Dolgov2012f}
 \textsc{S.\,V. Dolgov},  \textsc{B.\,N. Khoromskij},  and  \textsc{I.\,V.
  Oseledets},
Fast solution of parabolic problems in the tensor train/quantized tensor train
  format with initial application to the {Fokker}-{Planck} equation,
 \jr{SIAM J. Sci. Comput.} \textbf{34}(6), A3016--A3038 (2012).


\bibitem{Dolgov2012e}
 \textsc{S.\,V. Dolgov} and  \textsc{I.\,V. Oseledets},
Solution of linear systems and matrix inversion in the {TT}-format,
 \jr{SIAM J. Sci. Comput.} \textbf{34}(5), A2718--A2739 (2012).


\othercit
\bibitem{Dolgov2013a}
 \textsc{S.\,V. Dolgov} and  \textsc{D.\,V. Savostyanov},
Alternating minimal energy methods for linear systems in higher dimensions.
  {Part I}: {SPD} systems,
{arXiv} preprint 1301.6068, 2013.


\othercit
\bibitem{Dolgov2013}
 \textsc{S.\,V. Dolgov} and  \textsc{D.\,V. Savostyanov},
Alternating minimal residual methods for the solution of high-dimensional
  linear systems in the tensor train format,
in preparation, 2013.


\bibitem{Doostan2009}
 \textsc{A.~Doostan} and  \textsc{G.~Iaccarino},
A least-squares approximation of partial differential equations with
  high-dimensional random inputs,
 \jr{J. Comput. Phys.} \textbf{228}(12), 4332--4345 (2009).


\othercit
\bibitem{Doostan2012}
 \textsc{A.~Doostan},  \textsc{A.~Validi},  and  \textsc{G.~Iaccarino},
Non-intrusive low-rank separated approximation of high-dimensional stochastic
  models,
{arXiv} preprint 1210.1532, 2012.


\bibitem{Drineas2007}
 \textsc{P.~Drineas} and  \textsc{M.\,W. Mahoney},
A randomized algorithm for a tensor-based generalization of the singular value
  decomposition,
 \jr{Linear Algebra Appl.} \textbf{420}(2--3), 553--571 (2007).


\bibitem{Elden2009}
 \textsc{L.~Eld{\'e}n} and  \textsc{B.~Savas},
A {N}ewton-{G}rassmann method for computing the best multilinear rank-{$(r\sb
  1,r\sb 2,r\sb 3)$} approximation of a tensor,
 \jr{SIAM J. Matrix Anal. Appl.} \textbf{31}(2), 248--271 (2009).


\othercit
\bibitem{Espig2008}
 \textsc{M.~Espig},
Effiziente Bestapproximation mittels Summen von Elementartensoren in hohen
  Dimensionen,
PhD thesis, Fakult\"at f\"ur Mathematik und Informatik, Universit\"at Leipzig,
  2008.


\bibitem{Espig2009}
 \textsc{M.~Espig},  \textsc{L.~Grasedyck},  and  \textsc{W.~Hackbusch},
Black box low tensor-rank approximation using fiber-crosses,
 \jr{Constr. Appr.} \textbf{30}(3), 557--597 (2009).


\bibitem{Espig2012}
 \textsc{M.~Espig} and  \textsc{W.~Hackbusch},
A regularized {N}ewton method for the efficient approximation of tensors
  represented in the canonical tensor format,
 \jr{Numer. Math.} \textbf{122}(3), 489--525 (2012).


\bibitem{Espig2011}
 \textsc{M.~Espig},  \textsc{W.~Hackbusch},  \textsc{S.~Handschuh},  and
  \textsc{R.~Schneider},
Optimization problems in contracted tensor networks,
 \jr{Comput. Vis. Sci.} \textbf{14}, 271--285 (2011).


\othercit
\bibitem{Espig2013}
 \textsc{M.~Espig},  \textsc{W.~Hackbusch},  \textsc{A.~Litvinenko},
  \textsc{H.~Matthies},  and  \textsc{E.~Zander},
Efficient analysis of high dimensional data in tensor formats,
 in: Sparse Grids and Applications, edited by J.~Garcke and M.~Griebel, Lecture
  Notes in Computational Science and Engineering Vol.\,88 (Springer, 2013),
  pp.\,31--56.


\bibitem{Espig2013a}
 \textsc{M.~Espig},  \textsc{W.~Hackbusch},  \textsc{A.~Litvinenko},
  \textsc{H.\,G. Matthies},  and  \textsc{P.~W\"ahnert},
Efficient low-rank approximation of the stochastic {G}alerkin matrix in tensor
  formats,
 \jr{Comput. Math. Appl.} (2013),
To appear. Available at \url{http://dx.doi.org/10.1016/j.camwa.2012.10.008}.


\bibitem{Espig2012a}
 \textsc{M.~Espig},  \textsc{W.~Hackbusch},  \textsc{T.~Rohwedder},  and
  \textsc{R.~Schneider},
Variational calculus with sums of elementary tensors of fixed rank,
 \jr{Numer. Math.} \textbf{122}(52), 469--488 (2012).


\othercit
\bibitem{Espig2012b}
 \textsc{M.~Espig},  \textsc{M.~Schuster},  \textsc{A.~Killaitis},
  \textsc{N.~Waldren},  \textsc{P.~W\"ahnert},  \textsc{S.~Handschuh},  and
  \textsc{H.~Auer},
{TensorCalculus library}, 2012,
Available at \url{http://gitorious.org/tensorcalculus}.


\othercit
\bibitem{Exl2012}
 \textsc{L.~Exl},  \textsc{C.~Abert},  \textsc{N.\,J. Mauser},
  \textsc{T.~Schrefl},  \textsc{H.\,P. Stimming},  and  \textsc{D.~Suess},
{FFT-based Kronecker} product approximation to micromagnetic long-range
  interactions,
Tech. Rep. arXiv:1212.3509, 2012.


\bibitem{Exl2012a}
 \textsc{L.~Exl},  \textsc{W.~Auzinger},  \textsc{S.~Bance},
  \textsc{M.~Gusenbauer},  \textsc{F.~Reichel},  and  \textsc{T.~Schrefl},
Fast stray field computation on tensor grids,
 \jr{J. Comput. Phys.} \textbf{231}(7), 2840--2850 (2012).


\bibitem{Falco2012}
 \textsc{A.~Falc\'o} and  \textsc{W.~Hackbusch},
On minimal subspaces in tensor representations,
 \jr{Found. Comput. Math.} \textbf{12}, 765--803 (2012).


\othercit
\bibitem{Falco2013}
 \textsc{A.~Falc\'o},  \textsc{W.~Hackbusch},  and  \textsc{A.~Nouy},
Geometric structures in tensor representations,
Tech. Rep.~9, MPI MIS Leipzig, 2013.


\bibitem{Falco2011}
 \textsc{A.~Falc{\'o}} and  \textsc{A.~Nouy},
A proper generalized decomposition for the solution of elliptic problems in
  abstract form by using a functional {E}ckart-{Y}oung approach,
 \jr{J. Math. Anal. Appl.} \textbf{376}(2), 469--480 (2011).


\bibitem{Falco2012a}
 \textsc{A.~Falc{\'o}} and  \textsc{A.~Nouy},
Proper generalized decomposition for nonlinear convex problems in tensor
  {B}anach spaces,
 \jr{Numer. Math.} \textbf{121}, 503--530 (2012).


\bibitem{Figueroa2012}
 \textsc{L.\,E. Figueroa} and  \textsc{E.~S{\"u}li},
Greedy approximation of high-dimensional {O}rnstein-{U}hlenbeck operators,
 \jr{Found. Comput. Math.} \textbf{12}(5), 573--623 (2012).


\othercit
\bibitem{Flad2010}
 \textsc{H.\,J. Flad},  \textsc{W.~Hackbusch},  \textsc{B.\,N. Khoromskij},
  and  \textsc{R.~Schneider},
Concepts of data-sparse tensor-product approximation in many-particle
  modelling,
 in: Matrix methods: theory, algorithms and applications,  (World Sci. Publ.,
  Hackensack, NJ, 2010),  pp.\,313--347.


\bibitem{Flad2008}
 \textsc{H.\,J. Flad},  \textsc{B.\,N. Khoromskij},  \textsc{D.\,V.
  Savostyanov},  and  \textsc{E.\,E. Tyrtyshnikov},
Verification of the cross {3D} algorithm on quantum chemistry data,
 \jr{Rus. J. Numer. Anal. Math. Model.} \textbf{23}(4), 329--344 (2008).


\bibitem{Friedland2011}
 \textsc{S.~Friedland},  \textsc{V.~Mehrmann},  \textsc{A.~Miedlar},  and
  \textsc{M.~Nkengla},
Fast low rank approximations of matrices and tensors,
 \jr{Electron. J. Linear Algebra} \textbf{22}, 1031--1048 (2011).


\othercit
\bibitem{Friedland2012}
 \textsc{S.~Friedland},  \textsc{V.~Mehrmann},  \textsc{R.~Pajarola},  and
  \textsc{S.\,K. Suter},
On best rank one approximation of tensors,
{arXiv} preprint 1112.5914, 2012.


\bibitem{Frowis2010}
 \textsc{F.~Fr\"owis},  \textsc{V.~Nebendahl},  and  \textsc{W.~D\"ur},
Tensor operators: Constructions and applications for long-range interaction
  systems,
 \jr{Phys. Rev. A} \textbf{81}(Jun), 062337 (2010).


\bibitem{Gandy2011}
 \textsc{S.~Gandy},  \textsc{B.~Recht},  and  \textsc{I.~Yamada},
Tensor completion and low-n-rank tensor recovery via convex optimization,
 \jr{Inverse Problems} \textbf{27}(2), 025010 (2011).


\bibitem{Gavrilyuk2005}
 \textsc{I.\,P. Gavrilyuk},  \textsc{W.~Hackbusch},  and  \textsc{B.\,N.
  Khoromskij},
Hierarchical tensor-product approximation to the inverse and related operators
  for high-dimensional elliptic problems,
 \jr{Computing} \textbf{74}, 131--157 (2005).


\bibitem{Gavrilyuk2011}
 \textsc{I.\,P. Gavrilyuk} and  \textsc{B.\,N. Khoromskij},
Quantized-{TT}-{Cayley} transform to compute dynamics and spectrum of
  high-dimensional {Hamiltonians},
 \jr{Comput. Meth. Appl. Math.} \textbf{1}(3), 273--290 (2011).


\othercit
\bibitem{Giraldi2012}
 \textsc{L.~Giraldi},
Contributions aux M\'ethodes de Calcul Bas\'ees sur l'Approximation de Tenseurs
  et Applications en M\'ecanique Num\'erique,
PhD thesis, \'Ecole Centrale de Nantes, 2012,
Prelimary version.


\bibitem{Giraldi2013}
 \textsc{L.~Giraldi},  \textsc{A.~Nouy},  \textsc{G.~Legrain},  and
  \textsc{P.~Cartraud},
Tensor-based methods for numerical homogenization from high-resolution images,
 \jr{Comput. Meth. Appl. Mech. Engrg.} \textbf{254}, 154--169 (2013).


\othercit
\bibitem{Golub1996}
 \textsc{G.\,H. Golub} and  \textsc{C.\,F. Van~Loan},
Matrix Computations, third edition (Johns Hopkins University Press, Baltimore,
  MD, 1996).


\bibitem{Goreinov2012}
 \textsc{S.\,A. Goreinov},  \textsc{I.\,V. Oseledets},  and  \textsc{D.\,V.
  Savostyanov},
{Wedderburn} rank reduction and {Krylov} subspace method for tensor
  approximation. {Part} 1: {Tucker} case,
 \jr{SIAM J. Sci. Comput.} \textbf{34}(1), A1--A27 (2012).


\othercit
\bibitem{Goreinov2010}
 \textsc{S.\,A. Goreinov},  \textsc{I.\,V. Oseledets},  \textsc{D.\,V.
  Savostyanov},  \textsc{E.\,E. Tyrtyshnikov},  and  \textsc{N.\,L.
  Zamarashkin},
How to find a good submatrix,
 in: Matrix Methods: Theory, Algorithms, Applications,  (World Scientific,
  Hackensack, NY, 2010),  pp.\,247--256.


\bibitem{Goreinov1997}
 \textsc{S.\,A. Goreinov},  \textsc{E.\,E. Tyrtyshnikov},  and  \textsc{N.\,L.
  Zamarashkin},
A theory of pseudoskeleton approximations,
 \jr{Linear Algebra Appl.} \textbf{261}, 1--21 (1997).


\bibitem{Grasedyck2004}
 \textsc{L.~Grasedyck},
Existence and computation of low {K}ronecker-rank approximations for large
  linear systems of tensor product structure,
 \jr{Computing} \textbf{72}(3--4), 247--265 (2004).


\othercit
\bibitem{Grasedyck2010a}
 \textsc{L.~Grasedyck},
Hierarchical low rank approximation of tensors and multivariate functions,
  2010,
Lecture notes of Zurich summer school on "Sparse Tensor Discretizations of
  High-Dimensional Problems".


\bibitem{Grasedyck2010}
 \textsc{L.~Grasedyck},
Hierarchical singular value decomposition of tensors,
 \jr{SIAM J. Matrix Anal. Appl.} \textbf{31}(4), 2029--2054 (2010).


\othercit
\bibitem{Grasedyck2010b}
 \textsc{L.~Grasedyck},
Polynomial approximation in hierarchical {T}ucker format by
  vector-tensorization,
Tech. rep., RWTH Aachen, Germany, 2010.


\bibitem{Grasedyck2011}
 \textsc{L.~Grasedyck} and  \textsc{W.~Hackbusch},
An introduction to hierarchical {(H-)} rank and {TT}-rank of tensors with
  examples,
 \jr{Comput. Meth. Appl. Math.} \textbf{11}(3), 291--304 (2011).


\othercit
\bibitem{Griebel2011}
 \textsc{M.~Griebel} and  \textsc{H.~Harbrecht},
Approximation of two-variate functions: singular value decomposition versus
  regular sparse grids,
Preprint 1109, INS, Universit\"at Bonn, 2011.


\bibitem{Hackbusch2011}
 \textsc{W.~Hackbusch},
Tensorisation of vectors and their efficient convolution,
 \jr{Numer. Math.} \textbf{119}, 465--488 (2011).


\othercit
\bibitem{Hackbusch2012}
 \textsc{W.~Hackbusch},
Tensor Spaces and Numerical Tensor Calculus (Springer, 2012).


\bibitem{Hackbusch2012b}
 \textsc{W.~Hackbusch},
{$L^\infty$} estimation of tensor truncations,
 \jr{Numer. Math.} (2013),
To appear. Available at
  \url{http://www.mis.mpg.de/publications/preprints/2012/prepr2012-17.html}.


\bibitem{Hackbusch2006}
 \textsc{W.~Hackbusch} and  \textsc{B.\,N. Khoromskij},
Low-rank {K}ronecker-product approximation to multi-dimensional nonlocal
  operators. {I}. {S}eparable approximation of multi-variate functions,
 \jr{Computing} \textbf{76}(3-4), 177--202 (2006).


\bibitem{Hackbusch2006a}
 \textsc{W.~Hackbusch} and  \textsc{B.\,N. Khoromskij},
Low-rank {K}ronecker-product approximation to multi-dimensional nonlocal
  operators. {II}. {HKT} representation of certain operators,
 \jr{Computing} \textbf{76}(3-4), 203--225 (2006).


\bibitem{Hackbusch2007}
 \textsc{W.~Hackbusch} and  \textsc{B.\,N. Khoromskij},
Tensor-product approximation to operators and functions in high dimensions,
 \jr{J. Complexity} \textbf{23}(4-6), 697--714 (2007).


\bibitem{Hackbusch2008a}
 \textsc{W.~Hackbusch} and  \textsc{B.\,N. Khoromskij},
Tensor-product approximation to multidimensional integral operators and
  {G}reen's functions,
 \jr{SIAM J. Matrix Anal. Appl.} \textbf{30}(3), 1233--1253 (2008).


\bibitem{Hackbusch2012a}
 \textsc{W.~Hackbusch},  \textsc{B.\,N. Khoromskij},  \textsc{S.\,A. Sauter},
  and  \textsc{E.\,E. Tyrtyshnikov},
Use of tensor formats in elliptic eigenvalue problems,
 \jr{Numer. Linear Algebra Appl.} \textbf{19}(1), 133--151 (2012).


\bibitem{Hackbusch2005}
 \textsc{W.~Hackbusch},  \textsc{B.\,N. Khoromskij},  and  \textsc{E.\,E.
  Tyrtyshnikov},
Hierarchical {K}ronecker tensor-product approximations,
 \jr{J. Numer. Math.} \textbf{13}(2), 119--156 (2005).


\bibitem{Hackbusch2008}
 \textsc{W.~Hackbusch},  \textsc{B.\,N. Khoromskij},  and  \textsc{E.\,E.
  Tyrtyshnikov},
Approximate iterations for structured matrices,
 \jr{Numer. Math.} \textbf{109}(3), 365--383 (2008).


\bibitem{Hackbusch2009b}
 \textsc{W.~Hackbusch} and  \textsc{S.~K{\"u}hn},
A new scheme for the tensor representation,
 \jr{J. Fourier Anal. Appl.} \textbf{15}(5), 706--722 (2009).


\bibitem{Haegeman2011}
 \textsc{J.~Haegeman},  \textsc{I.~Cirac},  \textsc{T.~Osborne},
  \textsc{I.~Pi\'zorn},  \textsc{H.~Verschelde},  and  \textsc{F.~Verstraete},
Time-dependent variational principle for quantum lattices,
 \jr{Phys. Rev. Lett.} \textbf{107}(7), 070601 (2011).


\othercit
\bibitem{Haegeman2012}
 \textsc{J.~Haegeman},  \textsc{M.~Mari\"en},  \textsc{T.\,J. Osborne},  and
  \textsc{F.~Verstraete},
Geometry of matrix product states: metric, parallel transport and curvature,
{arXiv} preprint 1210.7710, 2012.


\bibitem{Halko2011}
 \textsc{N.~Halko},  \textsc{P.\,G. Martinsson},  and  \textsc{J.\,A. Tropp},
Finding structure with randomness: probabilistic algorithms for constructing
  approximate matrix decompositions,
 \jr{SIAM Rev.} \textbf{53}(2), 217--288 (2011).


\othercit
\bibitem{Handschuh2012}
 \textsc{S.~Handschuh},
Changing the topology of tensor networks,
Tech. Rep.~15, MPI MIS Leipzig, 2012.


\othercit
\bibitem{Hegland2011}
 \textsc{M.~Hegland} and  \textsc{J.~Garcke},
On the numerical solution of the chemical master equation with sums of rank one
  tensors,
 in: Proceedings of the 15th Biennial Computational Techniques and Applications
  Conference, CTAC-2010, edited by W.~McLean and A.\,J. Roberts, ANZIAM J.
  Vol.\,52 (August 2011),  pp.\,C628--C643.


\bibitem{Holtz2010}
 \textsc{S.~Holtz},  \textsc{T.~Rohwedder},  and  \textsc{R.~Schneider},
On manifolds of tensors of fixed {TT}-rank,
 \jr{Numer. Math.} \textbf{120}(4), 701--731 (2010).


\bibitem{Holtz2012}
 \textsc{S.~Holtz},  \textsc{T.~Rohwedder},  and  \textsc{R.~Schneider},
The alternating linear scheme for tensor optimization in the tensor train
  format,
 \jr{SIAM J. Sci. Comput.} \textbf{34}(2), A683--A713 (2012).


\othercit
\bibitem{Horn2013}
 \textsc{R.\,A. Horn} and  \textsc{C.\,R. Johnson},
Matrix analysis, second edition (Cambridge University Press, 2013).


\bibitem{Huckle2012a}
 \textsc{T.~Huckle} and  \textsc{K.~Waldherr},
Subspace iteration methods in terms of matrix product states,
 \jr{PAMM} \textbf{12}(1), 641--642 (2012).


\bibitem{Huckle2012}
 \textsc{T.~Huckle},  \textsc{K.~Waldherr},  and
  \textsc{T.~Schulte-Herbr{\"u}ggen},
Computations in quantum tensor networks,
 \jr{Linear Algebra Appl.} \textbf{438}(2), 750--781 (2013).


\bibitem{Huckle2013}
 \textsc{T.~Huckle},  \textsc{K.~Waldherr},  and
  \textsc{T.~Schulte-Herbr{\"u}ggen},
Exploiting matrix symmetries and physical symmetries in matrix product states
  and tensor trains,
 \jr{Linear and Multilinear Algebra} \textbf{61}(1), 91--122 (2013).


\bibitem{Ibragimov2009}
 \textsc{I.~Ibragimov} and  \textsc{S.~Rjasanow},
Three way decomposition for the {B}oltzmann equation,
 \jr{J. Comput. Math.} \textbf{27}(2-3), 184--195 (2009).


\othercit
\bibitem{Ishteva2011a}
 \textsc{M.~Ishteva},  \textsc{P.\,A. Absil},  and  \textsc{P.~Van~Dooren},
Jacobi algorithm for the best low multilinear rank approximation of symmetric
  tensors,
{Tech. rep. UCL-INMA-2011.011}, UC Louvain, Belgium, 2011.


\bibitem{Ishteva2011}
 \textsc{M.~Ishteva},  \textsc{P.\,A. Absil},  \textsc{S.~Van~Huffel},  and
  \textsc{L.~De~Lathauwer},
Best low multilinear rank approximation of higher-order tensors, based on the
  {R}iemannian trust-region scheme,
 \jr{SIAM J. Matrix Anal. Appl.} \textbf{32}(1), 115--135 (2011).


\bibitem{Ishteva2009}
 \textsc{M.~Ishteva},  \textsc{L.~De~Lathauwer},  \textsc{P.\,A. Absil},  and
  \textsc{S.~Van~Huffel},
Differential-geometric {N}ewton method for the best rank-{$(R_1,R_2,R_3)$}
  approximation of tensors,
 \jr{Numer. Algorithms} \textbf{51}(2), 179--194 (2009).


\bibitem{Jahnke2008}
 \textsc{T.~Jahnke} and  \textsc{W.~Huisinga},
A dynamical low-rank approach to the chemical master equation,
 \jr{Bulletin of Mathematical Biology} \textbf{70}, 2283--2302 (2008).


\othercit
\bibitem{Jondeau2010}
 \textsc{E.~Jondeau},  \textsc{E.~Jurczenko},  and  \textsc{M.~Rockinger},
Moment component analysis: An illustration with international stock markets,
{Swiss Finance Institute Research Paper No. 10-43}, 2010,
Available at \url{http://dx.doi.org/10.2139/ssrn.1694643}.


\othercit
\bibitem{Kazeev2013}
 \textsc{V.~Kazeev},  \textsc{M.~Khammash},  \textsc{M.~Nip},  and
  \textsc{C.~Schwab},
Direct solution of the chemical master equation using quantized tensor trains,
Tech. Rep. 2013-04, Seminar for Applied Mathematics, ETH Z{\"u}rich, 2013.


\othercit
\bibitem{Kazeev2011}
 \textsc{V.~Kazeev},  \textsc{B.\,N. Khoromskij},  and  \textsc{E.\,E.
  Tyrtyshnikov},
Multilevel {Toeplitz} matrices generated by {QTT} tensor-structured vectors and
  convolution with logarithmic complexity,
Tech. Rep.~36, MPI MIS Leipzig, 2011.


\othercit
\bibitem{Kazeev2012b}
 \textsc{V.~Kazeev},  \textsc{O.~Reichmann},  and  \textsc{C.~Schwab},
$hp$-{DG}-{QTT} solution of high-dimensional degenerate diffusion equations,
Tech. Rep. 2012-11, Seminar for Applied Mathematics, ETH Zurich, 2012.


\othercit
\bibitem{Kazeev2012a}
 \textsc{V.~Kazeev},  \textsc{O.~Reichmann},  and  \textsc{C.~Schwab},
Low-rank tensor structure of linear diffusion operators in the {TT} and {QTT}
  formats,
Tech. Rep. 2012-13, Seminar for Applied Mathematics, ETH Zurich, 2012.


\bibitem{Kazeev2012}
 \textsc{V.\,A. Kazeev} and  \textsc{B.\,N. Khoromskij},
Low-rank explicit {QTT} representation of {Laplace} operator and its inverse,
 \jr{SIAM J. Matrix Anal. Appl.} \textbf{33}(3), 742--758 (2012).


\bibitem{Kazeev2010}
 \textsc{V.\,A. Kazeev} and  \textsc{E.\,E. Tyrtyshnikov},
Structure of the {Hessian} matrix and an economical implementation of
  {Newton's} method in the problem of canonical approximation of tensors,
 \jr{Comput. Math. Math. Phys.} \textbf{50}(6), 927--945 (2010).


\bibitem{Khoromskaia2010}
 \textsc{V.~Khoromskaia},
Computation of the {H}artree-{F}ock exchange in the tensor-structured format.,
 \jr{Comput. Meth. Appl. Math.} \textbf{10}(2), 204--218 (2010).


\othercit
\bibitem{Khoromskaia2010a}
 \textsc{V.~Khoromskaia},
Numerical solution of the {Hartree-Fock} equation by multilevel
  tensor-structured methods,
PhD thesis, Technische Universit\"at Berlin, 2010.


\bibitem{Khoromskaia2012}
 \textsc{V.~Khoromskaia},  \textsc{D.~Andrae},  and  \textsc{B.\,N.
  Khoromskij},
Fast and accurate {3D} tensor calculation of the {Fock} operator in a general
  basis,
 \jr{Comput. Phys. Commun.} \textbf{183}(11), 2392--2404 (2012).


\bibitem{Khoromskaia2011}
 \textsc{V.~Khoromskaia},  \textsc{B.\,N. Khoromskij},  and
  \textsc{R.~Schneider},
{QTT} representation of the {Hartree} and exchange operators in electronic
  structure calculations,
 \jr{Comput. Meth. Appl. Math.} \textbf{11}(3), 327--341 (2011).


\othercit
\bibitem{Khoromskaia2013}
 \textsc{V.~Khoromskaia} and  \textsc{B.~Khoromskij},
{M{\o}ller-Plesset (MP2)} energy correction using tensor factorizations of the
  grid-based two-electron integrals,
Tech. Rep.~26, MPI MIS Leipzig, 2013.


\othercit
\bibitem{Khoromskaia2012c}
 \textsc{V.~Khoromskaia},  \textsc{B.~Khoromskij},  and  \textsc{R.~Schneider},
Tensor-structured factorized calculation of two-electron integrals in a general
  basis,
Tech. Rep.~29, MPI MIS Leipzig, 2012.


\othercit
\bibitem{Khoromskij2005}
 \textsc{B.\,N. Khoromskij},
An introduction to structured tensor-product approximation of discrete nonlocal
  operators, 2005,
Lecture notes No. 27, Leipzig.


\bibitem{Khoromskij2006}
 \textsc{B.\,N. Khoromskij},
Structured rank-$(r_1,\ldots,r_d)$ decomposition of function-related tensors in
  $\mathbb{R}^d$,
 \jr{Comput. Meth. Appl. Math.} \textbf{6}(2), 194--220 (2006).


\bibitem{Khoromskij2007}
 \textsc{B.\,N. Khoromskij},
Structured data-sparse approximation to high order tensors arising from the
  deterministic {Boltzmann} equation,
 \jr{Math. Comp.} \textbf{76}(259), 1291--1315 (2007).


\bibitem{Khoromskij2008}
 \textsc{B.\,N. Khoromskij},
On tensor approximation of {Green} iterations for {Kohn-Sham} equations,
 \jr{Comput. Vis. Sci.} \textbf{11}, 259--271 (2008).


\bibitem{Khoromskij2009b}
 \textsc{B.\,N. Khoromskij},
Tensor-structured preconditioners and approximate inverse of elliptic operators
  in $\mathbb{R}^d$,
 \jr{Constr. Approx.} \textbf{30}, 599--620 (2009).


\bibitem{Khoromskij2010b}
 \textsc{B.\,N. Khoromskij},
Fast and accurate tensor approximation of a multivariate convolution with
  linear scaling in dimension,
 \jr{J. Comput. Appl. Math.} \textbf{234}(11), 3122--3139 (2010).


\othercit
\bibitem{Khoromskij2011}
 \textsc{B.\,N. Khoromskij},
Introduction to tensor numerical methods in scientific computing, 2011,
Lecture notes, available at
  \url{http://www.math.uzh.ch/fileadmin/math/preprints/06_11.pdf}.


\bibitem{Khoromskij2011c}
 \textsc{B.\,N. Khoromskij},
{$O(d \log N)$}-quantics approximation of {$N-d$} tensors in high-dimensional
  numerical modeling,
 \jr{Constr. Approx.} \textbf{34}, 257--280 (2011).


\bibitem{Khoromskij2012}
 \textsc{B.\,N. Khoromskij},
Tensors-structured numerical methods in scientific computing: Survey on recent
  advances,
 \jr{Chemometrics and Intelligent Laboratory Systems} \textbf{110}(1), 1--19
  (2012).


\bibitem{Khoromskij2007a}
 \textsc{B.\,N. Khoromskij} and  \textsc{V.~Khoromskaia},
Low rank {T}ucker-type tensor approximation to classical potentials,
 \jr{Central European Journal of Mathematics} \textbf{5}, 523--550 (2007).


\bibitem{Khoromskij2009}
 \textsc{B.\,N. Khoromskij} and  \textsc{V.~Khoromskaia},
Multigrid accelerated tensor approximation of function related multidimensional
  arrays,
 \jr{SIAM J. Sci. Comput.} \textbf{31}(4), 3002--3026 (2009).


\bibitem{Khoromskij2009a}
 \textsc{B.\,N. Khoromskij},  \textsc{V.~Khoromskaia},  \textsc{S.\,R.
  Chinnamsetty},  and  \textsc{H.\,J. Flad},
Tensor decomposition in electronic structure calculations on {3D} cartesian
  grids,
 \jr{J. Comput. Phys.} \textbf{228}(16), 5749--5762 (2009).


\bibitem{Khoromskij2011d}
 \textsc{B.\,N. Khoromskij},  \textsc{V.~Khoromskaia},  and  \textsc{H.\,J.
  Flad},
Numerical solution of the {Hartree}--{Fock} equation in multilevel
  tensor-structured format,
 \jr{SIAM J. Sci. Comput.} \textbf{33}(1), 45--65 (2011).


\othercit
\bibitem{Khoromskij2010a}
 \textsc{B.\,N. Khoromskij} and  \textsc{I.\,V. Oseledets},
{DMRG+QTT} approach to computation of the ground state for the molecular
  {Schr\"odinger} operator,
Tech. Rep.~69, MPI MIS Leipzig, 2010.


\bibitem{Khoromskij2010}
 \textsc{B.\,N. Khoromskij} and  \textsc{I.\,V. Oseledets},
Quantics-{TT} collocation approximation of parameter-dependent and stochastic
  elliptic {PDEs},
 \jr{Comput. Meth. Appl. Math.} \textbf{10}(4), 376--394 (2010).


\bibitem{Khoromskij2011a}
 \textsc{B.\,N. Khoromskij} and  \textsc{I.\,V. Oseledets},
{QTT} approximation of elliptic solution operators in higher dimensions,
 \jr{Rus. J. Numer. Anal. Math. Model.} \textbf{26}(3), 303--322 (2011).


\othercit
\bibitem{Khoromskij2012a}
 \textsc{B.\,N. Khoromskij},  \textsc{I.\,V. Oseledets},  and
  \textsc{R.~Schneider},
Efficient time-stepping scheme for dynamics on {TT}-manifolds,
Tech. Rep.~24, MPI MIS Leipzig, 2012.


\bibitem{Khoromskij2011f}
 \textsc{B.\,N. Khoromskij},  \textsc{S.\,A. Sauter},  and  \textsc{A.~Veit},
Fast quadrature techniques for retarded potentials based on {TT}/{QTT} tensor
  approximation,
 \jr{Comput. Meth. Appl. Math.} \textbf{11}(3), 342--362 (2011).


\bibitem{Khoromskij2011b}
 \textsc{B.\,N. Khoromskij} and  \textsc{C.~Schwab},
Tensor-structured {G}alerkin approximation of parametric and stochastic
  elliptic {PDEs},
 \jr{SIAM J. Sci. Comput.} \textbf{33}(1), 364--385 (2011).


\bibitem{Koch2007}
 \textsc{O.~Koch} and  \textsc{C.~Lubich},
Dynamical low-rank approximation,
 \jr{SIAM J. Matrix Anal. Appl.} \textbf{29}(2), 434--454 (2007).


\bibitem{Koch2010}
 \textsc{O.~Koch} and  \textsc{C.~Lubich},
Dynamical tensor approximation,
 \jr{SIAM J. Matrix Anal. Appl.} \textbf{31}(5), 2360--2375 (2010).


\bibitem{Kolda2009}
 \textsc{T.\,G. Kolda} and  \textsc{B.\,W. Bader},
Tensor decompositions and applications,
 \jr{SIAM Review} \textbf{51}(3), 455--500 (2009).


\othercit
\bibitem{Kressner2012a}
 \textsc{D.~Kressner},  \textsc{M.~Ple\v{s}inger},  and  \textsc{C.~Tobler},
A preconditioned low-rank {CG} method for parameter-dependent {L}yapunov matrix
  equations,
Tech. rep., MATHICSE, EPF Lausanne, Switzerland, 2012,
Available at \url{http://anchp.epfl.ch}.


\bibitem{Kressner2010}
 \textsc{D.~Kressner} and  \textsc{C.~Tobler},
Krylov subspace methods for linear systems with tensor product structure,
 \jr{SIAM J. Matrix Anal. Appl.} \textbf{31}(4), 1688--1714 (2010).


\bibitem{Kressner2011}
 \textsc{D.~Kressner} and  \textsc{C.~Tobler},
Low-rank tensor {K}rylov subspace methods for parametrized linear systems,
 \jr{SIAM J. Matrix Anal. Appl.} \textbf{32}(4), 1288--1316 (2011).


\bibitem{Kressner2011a}
 \textsc{D.~Kressner} and  \textsc{C.~Tobler},
Preconditioned low-rank methods for high-dimensional elliptic {PDE} eigenvalue
  problems,
 \jr{Comput. Meth. Appl. Math.} \textbf{11}(3), 363--381 (2011).


\othercit
\bibitem{Kressner2012}
 \textsc{D.~Kressner} and  \textsc{C.~Tobler},
{\tt htucker} -- a {\sc matlab} toolbox for tensors in hierarchical {T}ucker
  format,
Tech. rep., MATHICSE, EPF Lausanne, Switzerland, 2012,
Available at \url{http://anchp.epfl.ch/htucker}.


\othercit
\bibitem{Kroonenberg2008}
 \textsc{P.\,M. Kroonenberg},
Applied Multiway Data Analysis (Wiley, 2008).


\othercit
\bibitem{Kuhn2012}
 \textsc{S.~K{\"u}hn},
{Hierarchische Tensordarstellung},
Dissertation, Universit\"at Leipzig, 2012.


\bibitem{Lamari2010}
 \textsc{H.~Lamari},  \textsc{A.~Ammar},  \textsc{P.~Cartraud},
  \textsc{G.~Legrain},  \textsc{F.~Chinesta},  and  \textsc{F.~Jacquemin},
Routes for efficient computational homogenization of nonlinear materials using
  the proper generalized decompositions,
 \jr{Arch. Comput. Methods Eng.} \textbf{17}, 373--391 (2010).


\othercit
\bibitem{Landsberg2012a}
 \textsc{J.\,M. Landsberg},
Tensors: geometry and applications, Graduate Studies in Mathematics,  Vol.\,128
  (American Mathematical Society, Providence, RI, 2012).


\bibitem{Landsberg2012}
 \textsc{J.\,M. Landsberg},  \textsc{Y.~Qi},  and  \textsc{K.~Ye},
On the geometry of tensor network states,
 \jr{Quantum Inf. Comput.} \textbf{12}(3-4), 346--354 (2012).


\bibitem{LeBris2009}
 \textsc{C.~Le~Bris},  \textsc{T.~Leli{\`e}vre},  and  \textsc{Y.~Maday},
Results and questions on a nonlinear approximation approach for solving
  high-dimensional partial differential equations,
 \jr{Constr. Approx.} \textbf{30}(3), 621--651 (2009).


\bibitem{Lebedeva2010}
 \textsc{O.\,S. Lebedeva},
Block tensor conjugate gradient-type method for {Rayleigh} quotient
  minimization in two-dimensional case,
 \jr{Comput. Math. Math. Phys.} \textbf{50}(5), 749--765 (2010).


\bibitem{Lebedeva2011}
 \textsc{O.\,S. Lebedeva},
Tensor conjugate-gradient-type method for {R}ayleigh quotient minimization in
  block {QTT}-format,
 \jr{Russian J. Numer. Anal. Math. Modelling} \textbf{26}(5), 465--489 (2011).


\othercit
\bibitem{Liu2009}
 \textsc{J.~Liu},  \textsc{P.~Musialski},  \textsc{P.~Wonka},  and
  \textsc{J.~Ye},
Tensor completion for estimating missing values in visual data,
 in: Computer Vision, 2009 IEEE 12th International Conference on,  (29
  2009-oct. 2 2009),  pp.\,2114--2121.


\bibitem{Lubich2005}
 \textsc{C.~Lubich},
On variational approximations in quantum molecular dynamics,
 \jr{Math. Comp.} \textbf{74}(250), 765--779 (2005).


\othercit
\bibitem{Lubich2008}
 \textsc{C.~Lubich},
From quantum to classical molecular dynamics: reduced models and numerical
  analysis, Zurich Lectures in Advanced Mathematics (European Mathematical
  Society (EMS), Z\"urich, 2008).


\othercit
\bibitem{Lubich2013}
 \textsc{C.~Lubich} and  \textsc{I.\,V. Oseledets},
A projector-splitting integrator for dynamical low-rank approximation,
{arXiv} preprint 1301.5512, 2013.


\othercit
\bibitem{Lubich2012}
 \textsc{C.~Lubich},  \textsc{T.~Rohwedder},  \textsc{R.~Schneider},  and
  \textsc{B.~Vandereycken},
Dynamical approximation of hierarchical {T}ucker and tensor-train tensors,
Tech. rep., July 2012.


\othercit
\bibitem{Mach2011}
 \textsc{T.~Mach},
Computing inner eigenvalues of matrices in tensor train matrix format,
Tech. Rep. MPIMD/11-09, MPI Magdeburg, 2011,
Available at \url{http://www.mpi-magdeburg.mpg.de/preprints/2011/09/}.


\othercit
\bibitem{Mach2012}
 \textsc{T.~Mach} and  \textsc{J.~Saak},
Towards an {ADI} iteration for tensor structured equations,
Tech. rep., MPI Magdeburg, 2012,
Available at \url{www.mpi-magdeburg.mpg.de/preprints/2011/12/}.


\bibitem{Matthies2012}
 \textsc{H.\,G. Matthies} and  \textsc{E.~Zander},
Solving stochastic systems with low-rank tensor compression,
 \jr{Linear Algebra Appl.} \textbf{436}(10), 3819--3838 (2012).


\othercit
\bibitem{Meszmer2012}
 \textsc{P.~Meszmer} and  \textsc{J.~Ballani},
{Tensor structured evaluation of singular volume integrals},
Tech. Rep.~70, {MPI MIS Leipzig}, 2012.


\bibitem{Meyer2011}
 \textsc{G.~Meyer},  \textsc{S.~Bonnabel},  and  \textsc{R.~Sepulchre},
Regression on fixed-rank positive semidefinite matrices: a {R}iemannian
  approach,
 \jr{J. Mach. Learn. Res.} \textbf{12}, 593--625 (2011).


\othercit
\bibitem{Meyer2009}
 \textsc{H.\,D. Meyer},  \textsc{F.~Gatti},  and  \textsc{G.\,A. Worth} (eds.),
Multidimensional Quantum Dynamics: MCTDH Theory and Applications (Wiley-VCH
  Verlag GmbH \& Co. KGaA, 2009).


\bibitem{Meyer1990}
 \textsc{H.\,D. Meyer},  \textsc{U.~Manthe},  and  \textsc{L.~Cederbaum},
The multi-configurational time-dependent hartree approach,
 \jr{Chemical Physics Letters} \textbf{165}(1), 73 -- 78 (1990).


\bibitem{Mohlenkamp2013}
 \textsc{M.\,J. Mohlenkamp},
Musings on multilinear fitting,
 \jr{Linear Algebra Appl.} \textbf{438}(2), 834--852 (2013).


\bibitem{Mohlenkamp2005}
 \textsc{M.\,J. Mohlenkamp} and  \textsc{L.~Monz{\'o}n},
Trigonometric identities and sums of separable functions,
 \jr{Math. Intelligencer} \textbf{27}(2), 65--69 (2005).


\bibitem{Mokdad2010}
 \textsc{B.~Mokdad},  \textsc{A.~Ammar},  \textsc{M.~Normandin},
  \textsc{F.~Chinesta},  and  \textsc{J.\,R. Clermont},
A fully deterministic micro-macro simulation of complex flows involving
  reversible network fluid models,
 \jr{Math. Comput. Simulation} \textbf{80}(9), 1936--1961 (2010).


\bibitem{Naraparaju2011}
 \textsc{K.\,K. Naraparaju} and  \textsc{J.~Schneider},
Generalized cross approximation for 3d-tensors,
 \jr{Comput. Vis. Sci.} \textbf{14}(3), 105--115 (2011).


\othercit
\bibitem{Ngo2012}
 \textsc{T.~Ngo} and  \textsc{Y.~Saad},
Scaled gradients on {G}rassmann manifolds for matrix completion,
{Preprint ys-2012-5}, Dept. Computer Science and Engineering, University of
  Minnesota, Minneapolis, MN, 2012.


\bibitem{Nonnenmacher2008}
 \textsc{A.~Nonnenmacher} and  \textsc{C.~Lubich},
Dynamical low-rank approximation: applications and numerical experiments,
 \jr{Mathematics and Computers in Simulation} \textbf{79}(4), 1346--1357
  (2008).


\bibitem{Nouy2007}
 \textsc{A.~Nouy},
A generalized spectral decomposition technique to solve a class of linear
  stochastic partial differential equations,
 \jr{Comput. Meth. Appl. Mech. Engrg.} \textbf{196}(45--48), 4521--4537 (2007).


\bibitem{Nouy2008}
 \textsc{A.~Nouy},
Generalized spectral decomposition method for solving stochastic finite element
  equations: invariant subspace problem and dedicated algorithms,
 \jr{Comput. Meth. Appl. Mech. Engrg.} \textbf{197}(51--52), 4718--4736 (2008).


\bibitem{Nouy2010}
 \textsc{A.~Nouy},
Proper generalized decompositions and separated representations for the
  numerical solution of high dimensional stochastic problems,
 \jr{Arch. Comput. Methods Eng.} \textbf{17}, 403--434 (2010).


\bibitem{Olshevsky2006}
 \textsc{V.~Olshevsky},  \textsc{I.\,V. Oseledets},  and  \textsc{E.\,E.
  Tyrtyshnikov},
Tensor properties of multilevel {Toeplitz} and related matrices,
 \jr{Linear Algebra Appl.} \textbf{412}(1), 1--21 (2006).


\othercit
\bibitem{Oseledets2013}
 \textsc{I.~Oseledets},  \textsc{S.~Dolgov},  and  \textsc{D.~Savostyanov},
{ttpy}, 2013,
Available at \url{https://github.com/oseledets/ttpy/}.


\bibitem{Oseledets2010b}
 \textsc{I.\,V. Oseledets},
Approximation of $2^d \times 2^d$ matrices using tensor decomposition,
 \jr{SIAM J. Matrix Anal. Appl.} \textbf{31}(4), 2130--2145 (2010).


\bibitem{Oseledets2011c}
 \textsc{I.\,V. Oseledets},
{DMRG} approach to fast linear algebra in the {TT}--format,
 \jr{Comput. Meth. Appl. Math} \textbf{11}(3), 382--393 (2011).


\othercit
\bibitem{Oseledets2011a}
 \textsc{I.\,V. Oseledets},
{MATLAB TT-Toolbox Version 2.2}, May 2011,
Available at \url{http://spring.inm.ras.ru/osel/?page_id=24}.


\bibitem{Oseledets2012a}
 \textsc{I.\,V. Oseledets},
Constructive representation of functions in low-rank tensor formats,
 \jr{Constr. Appr.} \textbf{37}(1), 1--18 (2013).


\bibitem{Oseledets2008}
 \textsc{I.\,V. Oseledets},  \textsc{D.\,V. Savostianov},  and  \textsc{E.\,E.
  Tyrtyshnikov},
{T}ucker dimensionality reduction of three-dimensional arrays in linear time,
 \jr{SIAM J. Matrix Anal. Appl.} \textbf{30}(3), 939--956 (2008).


\bibitem{Oseledets2009b}
 \textsc{I.\,V. Oseledets},  \textsc{D.\,V. Savostyanov},  and  \textsc{E.\,E.
  Tyrtyshnikov},
Linear algebra for tensor problems,
 \jr{Computing} \textbf{85}(3), 169--188 (2009).


\bibitem{Oseledets2010c}
 \textsc{I.\,V. Oseledets},  \textsc{D.\,V. Savostyanov},  and  \textsc{E.\,E.
  Tyrtyshnikov},
Cross approximation in tensor electron density computations,
 \jr{Numer. Linear Algebra Appl.} \textbf{17}(6), 935--952 (2010).


\bibitem{Oseledets2009a}
 \textsc{I.\,V. Oseledets} and  \textsc{E.\,E. Tyrtyshnikov},
Breaking the curse of dimensionality, or how to use {SVD} in many dimensions,
 \jr{SIAM J. Sci. Comput.} \textbf{31}(5), 3744--3759 (2009).


\othercit
\bibitem{Oseledets2009f}
 \textsc{I.\,V. Oseledets} and  \textsc{E.\,E. Tyrtyshnikov},
Tensor tree decomposition does not need a tree,
Preprint (Submitted to Linear Algebra Appl) 2009-04, INM RAS, Moscow, 2009.


\bibitem{Oseledets2010a}
 \textsc{I.\,V. Oseledets} and  \textsc{E.\,E. Tyrtyshnikov},
{TT}-cross approximation for multidimensional arrays,
 \jr{Linear Algebra Appl.} \textbf{432}(1), 70--88 (2010).


\bibitem{Oseledets2011b}
 \textsc{I.\,V. Oseledets} and  \textsc{E.\,E. Tyrtyshnikov},
Algebraic wavelet transform via quantics tensor train decomposition,
 \jr{SIAM J. Sci. Comput.} \textbf{33}(3), 1315--1328 (2011).


\bibitem{Oseledets2011d}
 \textsc{I.\,V. Oseledets},  \textsc{E.\,E. Tyrtyshnikov},  and  \textsc{N.\,L.
  Zamarashkin},
Tensor-train ranks of matrices and their inverses,
 \jr{Comput. Meth. Appl. Math} \textbf{11}(3), 394--403 (2011).


\bibitem{Oseledets2011}
 \textsc{I.~Oseledets},
Tensor-train decomposition,
 \jr{SIAM J. Sci. Comput.} \textbf{33}(5), 2295--2317 (2011).


\bibitem{Ostlund1995}
 \textsc{S.~\"Ostlund} and  \textsc{S.~Rommer},
Thermodynamic limit of density matrix renormalization,
 \jr{Phys. Rev. Lett.} \textbf{75}(Nov), 3537--3540 (1995).


\othercit
\bibitem{Pan2012}
 \textsc{V.\,Y. Pan},  \textsc{G.~Qian},  and  \textsc{A.\,L. Zheng},
Randomized matrix computations,
{arXiv} preprint 1210.7476, 2012.


\othercit
\bibitem{Penrose1971}
 \textsc{R.~Penrose},
Applications of negative dimensional tensors,
 in: Combinatorial {M}athematics and its {A}pplications ({P}roc. {C}onf.,
  {O}xford, 1969),  (Academic Press, London, 1971),  pp.\,221--244.


\bibitem{Perez-Garcia2007}
 \textsc{D.~Perez-Garcia},  \textsc{F.~Verstraete},  \textsc{M.\,M. Wolf},  and
   \textsc{J.\,I. Cirac},
Matrix product state representations,
 \jr{Quantum Info. Comput.} \textbf{7}(5), 401--430 (2007).


\othercit
\bibitem{Phan2012}
 \textsc{A.\,H. Phan},  \textsc{P.~Tichavsk\'y},  and  \textsc{A.~Cichocki},
{CANDECOMP/PARAFAC} decomposition of high-order tensors through tensor
  reshaping,
{arXiv} preprint 1211.3796, 2012.


\othercit
\bibitem{Phan2012b}
 \textsc{A.\,H. Phan},  \textsc{P.~Tichavsk\'y},  and  \textsc{A.~Cichocki},
On fast computation of gradients for {CANDECOMP/PARAFAC} algorithms,
{arXiv} preprint 1204.1586, 2012.


\bibitem{Phan2013}
 \textsc{A.\,H. Phan},  \textsc{P.~Tichavsk\'y},  and  \textsc{A.~Cichocki},
Low complexity damped gauss--newton algorithms for {CANDECOMP/PARAFAC},
 \jr{SIAM Journal on Matrix Analysis and Applications} \textbf{34}(1), 126--147
  (2013).


\bibitem{Rajih2008}
 \textsc{M.~Rajih},  \textsc{P.~Comon},  and  \textsc{R.\,A. Harshman},
Enhanced line search: a novel method to accelerate {PARAFAC},
 \jr{SIAM J. Matrix Anal. Appl.} \textbf{30}(3), 1128--1147 (2008).


\othercit
\bibitem{Rohwedder2011}
 \textsc{T.~Rohwedder} and  \textsc{A.~Uschmajew},
Local convergence of alternating schemes for optimization of convex problems in
  the {TT} format,
Tech. rep., 2011.


\bibitem{Sanz2009}
 \textsc{M.~Sanz},  \textsc{M.\,M. Wolf},  \textsc{D.~P\'erez-Garc\'ia},  and
  \textsc{J.\,I. Cirac},
Matrix product states: Symmetries and two-body {H}amiltonians,
 \jr{Phys. Rev. A} \textbf{79}(Apr), 042308 (2009).


\bibitem{Savas2013}
 \textsc{B.~Savas} and  \textsc{L.~Eld\'en},
Krylov-type methods for tensor computations {I},
 \jr{Linear Algebra Appl.} \textbf{438}(2), 891--918 (2013).


\bibitem{Savas2010}
 \textsc{B.~Savas} and  \textsc{L.\,H. Lim},
Quasi-{N}ewton methods on {G}rassmannians and multilinear approximations of
  tensors,
 \jr{SIAM J. Sci. Comput.} \textbf{32}(6), 3352--3393 (2010).


\bibitem{Savostyanov2012a}
 \textsc{D.\,V. Savostyanov},
{QTT}-rank-one vectors with {QTT}-rank-one and full-rank {Fourier} images,
 \jr{Linear Algebra Appl.} \textbf{436}(9), 3215--3224 (2012).


\othercit
\bibitem{Savostyanov2011a}
 \textsc{D.\,V. Savostyanov} and  \textsc{I.\,V. Oseledets},
Fast adaptive interpolation of multi-dimensional arrays in tensor train format,
 in: Proceedings of 7th International Workshop on Multidimensional Systems
  (nDS),  (IEEE, 2011).


\bibitem{Savostyanov2012}
 \textsc{D.\,V. Savostyanov},  \textsc{E.\,E. Tyrtyshnikov},  and
  \textsc{N.\,L. Zamarashkin},
Fast truncation of mode ranks for bilinear tensor operations,
 \jr{Numer. Linear Algebra Appl.} \textbf{19}(1), 103--111 (2012).


\othercit
\bibitem{Schneider2012}
 \textsc{R.~Schneider},
Tensor product approximation and numerical solution of the electronic
  {Schr\"odinger} equation, 2012,
John von Neumann Gastvorlesung, available at
  \url{http://www-m2.ma.tum.de/bin/view/Allgemeines/JohnVonNeumannLectureSS12}.


\othercit
\bibitem{Schneider2013}
 \textsc{R.~Schneider} and  \textsc{A.~Uschmajew},
Approximation rates for the hierarchical tensor format in certain periodic
  {S}obolev spaces,
Tech. rep., MATHICSE, EPFL, 2013.


\bibitem{Schollwock2011}
 \textsc{U.~Schollw\"ock},
The density-matrix renormalization group in the age of matrix product states,
 \jr{Annals of Physics} \textbf{326} (2011).


\bibitem{Sharma2012}
 \textsc{S.~Sharma} and  \textsc{G.\,K.\,L. Chan},
Spin-adapted density matrix renormalization group algorithms for quantum
  chemistry,
 \jr{J. Chem. Phys.} \textbf{136}(12), 124121 (2012).


\othercit
\bibitem{Sharma2013}
 \textsc{S.~Sharma} and  \textsc{G.\,K.\,L. Chan},
Block v 0.9.6, 2013,
Available at \url{http://www.princeton.edu/chemistry/chan/software/dmrg/}.


\bibitem{Shi2006}
 \textsc{Y.\,Y. Shi},  \textsc{L.\,M. Duan},  and  \textsc{G.~Vidal},
Classical simulation of quantum many-body systems with a tree tensor network,
 \jr{Physical Review A} \textbf{74}(2), 022320 (2006).


\othercit
\bibitem{Smilde2004}
 \textsc{A.~Smilde},  \textsc{R.~Bro},  and  \textsc{P.~Geladi},
Multi-way Analysis: Applications in the Chemical Sciences (Wiley, 2004).


\othercit
\bibitem{Sorber2013}
 \textsc{L.~Sorber},  \textsc{M.~Van~Barel},  and  \textsc{L.~De~Lathauwer},
Tensorlab v1.0, 2013,
Available online, \url{http://esat.kuleuven.be/sista/tensorlab/}.


\bibitem{Stegeman2010}
 \textsc{A.~Stegeman} and  \textsc{P.~Comon},
Subtracting a best rank-1 approximation may increase tensor rank,
 \jr{Linear Algebra Appl.} \textbf{433}(7), 1276--1300 (2010).


\bibitem{Takaaki2012}
 \textsc{N.~Takaaki},
Algebraic reconstruction of the general-order poles of a meromorphic function,
 \jr{Inverse Problems} \textbf{28}(2), 025008 (2012).


\bibitem{Temlyakov1986}
 \textsc{V.\,N. Temlyakov},
Approximations of functions with bounded mixed derivative,
 \jr{Trudy Mat. Inst. Steklov.} \textbf{178}, 113 (1986),
In Russian; English translation in Proc. Steklov Inst. Math. 1989, no. 1 (178).


\bibitem{Temlyakov1992}
 \textsc{V.\,N. Temlyakov},
Bilinear approximation and related questions,
 \jr{Trudy Mat. Inst. Steklov} pp.\,229--248 (1992),
In Russian; English translation: Proc. Steklov Inst. Math. 1993, no. 4 (194),
  245--265.


\bibitem{Temlyakov1992a}
 \textsc{V.\,N. Temlyakov},
Estimates for the best bilinear approximations of functions and approximation
  numbers of integral operators,
 \jr{Mat. Zametki} \textbf{51}(5), 125--134, 159 (1992),
In Russian; English translation in Math. Notes 1992, 51(5), 510--517.


\othercit
\bibitem{Tobler2012}
 \textsc{C.~Tobler},
Low-rank Tensor Methods for Linear Systems and Eigenvalue Problems,
PhD thesis, ETH Zurich, Switzerland, 2012.


\othercit
\bibitem{Tsourakakis2010}
 \textsc{C.\,E. Tsourakakis},
{MACH:} fast randomized tensor decompositions,
 in: SIAM International Conference on Data Mining,  (2010).


\bibitem{Tyrtyshnikov2000}
 \textsc{E.\,E. Tyrtyshnikov},
Incomplete cross approximation in the mosaic-skeleton method,
 \jr{Computing} \textbf{64}(4), 367--380 (2000),
International GAMM-Workshop on Multigrid Methods (Bonn, 1998).


\bibitem{Tyrtyshnikov2003}
 \textsc{E.\,E. Tyrtyshnikov},
Tensor approximations of matrices generated by asymptotically smooth functions,
 \jr{Mat. Sb.} \textbf{194}(6), 147--160 (2003).


\bibitem{Tyrtyshnikov2004}
 \textsc{E.\,E. Tyrtyshnikov},
Kronecker-product approximations for some function-related matrices,
 \jr{Linear Algebra Appl.} \textbf{379}, 423--437 (2004),
Tenth Conference of the International Linear Algebra Society.


\bibitem{Ullmann2010}
 \textsc{E.~Ullmann},
A {K}ronecker product preconditioner for stochastic {G}alerkin finite element
  discretizations,
 \jr{SIAM J. Sci. Comput.} \textbf{32}(2), 923--946 (2010).


\bibitem{Uschmajew2011}
 \textsc{A.~Uschmajew},
Regularity of tensor product approximations to square integrable functions,
 \jr{Constr. Approx.} \textbf{34}, 371--391 (2011).


\bibitem{Uschmajew2012a}
 \textsc{A.~Uschmajew},
Local convergence of the alternating least squares algorithm for canonical
  tensor approximation,
 \jr{SIAM J. Matrix Anal. Appl.} \textbf{33}(2), 639--652 (2012).


\othercit
\bibitem{Uschmajew2013}
 \textsc{A.~Uschmajew},
{Zur Theorie der Niedrigrangapproximation in Tensorprodukten von
  Hilbertr\"aumen},
PhD thesis, Technische Universit\"at Berlin, 2013.


\othercit
\bibitem{Uschmajew2012}
 \textsc{A.~Uschmajew} and  \textsc{B.~Vandereycken},
The geometry of algorithms using hierarchical tensors,
Tech. rep., January 2012.


\othercit
\bibitem{VanLoan1993}
 \textsc{C.\,F. Van~Loan} and  \textsc{N.~Pitsianis},
Approximation with {K}ronecker products,
 in: Linear algebra for large scale and real-time applications ({L}euven,
  1992), , NATO Adv. Sci. Inst. Ser. E Appl. Sci.,  Vol.\,232 (Kluwer Acad.
  Publ., Dordrecht, 1993),  pp.\,293--314.


\othercit
\bibitem{Vandereycken2012}
 \textsc{B.~Vandereycken},
Low-rank matrix completion by {R}iemannian optimization,
Tech. rep., MATHICSE, EPF Lausanne, Switzerland, 2012.


\bibitem{Vandereycken2010}
 \textsc{B.~Vandereycken} and  \textsc{S.~Vandewalle},
A {R}iemannian optimization approach for computing low-rank solutions of
  {L}yapunov equations,
 \jr{SIAM J. Matrix Anal. Appl.} \textbf{31}(5), 2553--2579 (2010).


\bibitem{Vannieuwenhoven2012}
 \textsc{N.~Vannieuwenhoven},  \textsc{R.~Vandebril},  and
  \textsc{K.~Meerbergen},
A new truncation strategy for the higher-order singular value decomposition,
 \jr{SIAM J. Sci. Comput.} \textbf{34}(2), A1027--A1052 (2012).


\othercit
\bibitem{Verstraete2004b}
 \textsc{F.~Verstraete} and  \textsc{J.\,I. Cirac},
Renormalization algorithms for quantum-many body systems in two and higher
  dimensions,
Tech. Rep. arXiv:cond-mat/0407066, 2004.


\bibitem{Verstraete2004a}
 \textsc{F.~Verstraete} and  \textsc{J.\,I. Cirac},
Valence-bond states for quantum computation,
 \jr{Phys. Rev. A} \textbf{70}(Dec), 060302 (2004).


\bibitem{Verstraete2004}
 \textsc{F.~Verstraete},  \textsc{J.\,J. Garc\'ia-Ripoll},  and  \textsc{J.\,I.
  Cirac},
Matrix product density operators: Simulation of finite-temperature and
  dissipative systems,
 \jr{Phys. Rev. Lett.} \textbf{93}(Nov), 207204 (2004).


\bibitem{Vidal2003}
 \textsc{G.~Vidal},
{Efficient Classical Simulation of Slightly Entangled Quantum Computations},
 \jr{Phys. Rev. Lett.} \textbf{91}, 147902 (2003).


\bibitem{Vidal2007}
 \textsc{G.~Vidal},
Entanglement renormalization,
 \jr{Phys. Rev. Lett.} \textbf{99}(Nov), 220405 (2007).


\bibitem{Wang2003}
 \textsc{H.~Wang} and  \textsc{M.~Thoss},
Multilayer formulation of the multiconfiguration time-dependent {H}artree
  theory,
 \jr{J. Chem. Phys} \textbf{119}(3), 1289--1299 (2003).


\bibitem{White1992}
 \textsc{S.\,R. White},
Density matrix formulation for quantum renormalization groups,
 \jr{Phys. Rev. Lett.} \textbf{69}, 2863--2866 (1992).


\othercit
\bibitem{Wise2013}
 \textsc{B.\,M. Wise} and  \textsc{N.\,B. Gallagher},
{PLS Toolbox 7.0.3}, 2013,
Available at \url{http://www.eigenvector.com}.


\bibitem{Worth2003}
 \textsc{G.\,A. Worth},  \textsc{M.\,H. Beck},  \textsc{A.~J\"ackle},  and
  \textsc{H.\,D. Meyer},
The {MCTDH} package, version 8.3 (2003),
See \url{http://www.pci.uni-heidelberg.de/tc/usr/mctdh/}.


\othercit
\bibitem{Zhou2012}
 \textsc{G.~Zhou},  \textsc{A.~Cichocki},  and  \textsc{S.~Xie},
Accelerated canonical polyadic decomposition by using mode reduction,
{arXiv} preprint 1211.3500, 2012.


\bibitem{Zwolak2004}
 \textsc{M.~Zwolak} and  \textsc{G.~Vidal},
Mixed-state dynamics in one-dimensional quantum lattice systems: A
  time-dependent superoperator renormalization algorithm,
 \jr{Phys. Rev. Lett.} \textbf{93}(Nov), 207205 (2004).


\end{thebibliography}

\end{document}